\documentstyle[12pt]{article} \textwidth485pt \hoffset-52pt
\def\N{{\bf N}} \def\Z{{\bf Z}} \def\R{{\bf R}} \def\C{{\bf C}}

\begin{document}
\baselineskip=0.5cm

$\:$ \smallskip
\begin{center}
 ON REMOVABLE SINGULARITIES\\ FOR CR FUNCTIONS \\ IN HIGHER
CODIMENSION
\end{center}

\begin{center}
{\bf J. Merker}
\end{center}

\bigskip
\bigskip
\bigskip
\smallskip
\smallskip

In recent years, several papers (for a complete reference list, see
Chirka and Stout \cite{CS}) have been published on the subject of
removable singularities for the boundary values of holomorphic
functions on some domains or hypersurfaces in the complex euclidean
space.

In this paper, we study the higher codimensional case. Our results for
the hypersurface case are weaker than those in \cite{CS} and
\cite{JO2}, for the smoothness assumption.

Note with $T_{z_0} M$ the usual tangent space of a real manifold
$M\subset \C^n$ at $z_0 \in M$ and by $T_{z_0}^c M= T_{z_0}M \cap
JT_{z_0}M$ its complex tangent space, where $J$ denotes the complex
structure on $T\C^n$. $M$ is said to be generic if $T_z M + J T_z M=
T_z\C^n$ for all $z \in M$. We consider continuous distributional
solutions of the tangential Cauchy-Riemann equations on $M$, which
will be referred as {\it CR functions} on $M$. Let $z_0 \in M$. By a
wedge of edge $M$ at $z_0$, we mean an open set in $\C^n$ of the form
$${\cal W}=\{z+\eta; \ z \in U, \eta \in C\},$$ for some open
neighborhood $U$ of $z_0$ in $M$ and some convex truncated open cone
$C$ in $T_{z_0}\C^n / T_{z_0}M$, {\em i.e.} the intersection of a
convex open cone with a ball centered at $0$.

Let now $N\ni z_0$ be a proper $C^1$ submanifold of $M$ and assume
that $M$ is minimal at $z_0$. Our submanifolds will always be assumed
to be {\it embedded submanifolds}. Call $N$ {\it removable at $z_0$}
if there exists a wedge ${\cal W}$ of edge $M$ at $z_0$ with CR
functions on $M\backslash N$ extending holomorphically into ${\cal W}$
and continuously in ${\cal W}\cup (M\backslash N)$.

Our main results are the following.

\smallskip
\noindent
{\sc Theorem 1.} {\it Let $M$ be a $C^{2,\alpha}$-smooth
$(0<\alpha<1)$ generic manifold in $\C^n (n\geq 2)$ minimal at every
point, with $\hbox{CRdim} \ M = p \geq 1$. Then every $C^1$
submanifold $N \subset M$ with $codim_M N \geq 3$ and $T_{z}N \not
\supset T_{z}^cM$ for every $z\in N$ is removable.}

\smallskip
\noindent
{\sc Theorem 2.} {\it Let $M$ be a $C^{2,\alpha}$-smooth
$(0<\alpha<1)$ generic manifold in $\C^n (n\geq 3)$ minimal at every
point, with $\hbox{CRdim} \ M = p \geq 2$. Then every connected $C^1$
submanifold $N\subset M$ with $codim_M N=2$ and $T_{z}N\not\supset
T_{z}^cM$ for every $z \in N$ is removable provided $N$ does not
consist of a CR manifold with $\hbox{CRdim} \ N=p-1$.}

\smallskip
Minimality is understood in the sense of Tumanov. The following
theorem is due to J\"oricke, via a minimalization theorem
\cite{JO3}. We recall her proof and give another one.

\smallskip
\noindent
{\sc Theorem 3.} {\it Let $M$ be as in Theorem 2 and let $N$ be a
connected $C^2$ submanifold with $\hbox{codim}_M N =1$ which is
generic in $\C^n$. Let $z_0\in N$, let $\Phi $ be a closed subset of
$N$ with $\Phi \neq N$ and $b\Phi \ni z_0$. If $\Phi$ does not contain
the germ through $z_0$ of a CR manifold $\Sigma$ with $\hbox{CRdim} \
\Sigma= p-1$, $\Phi$ is removable at $z_0$.}

\smallskip
It should be noted that, though not being a CR manifold, a
two-codimensional manifold $N$ can contain proper submanifolds
$\Sigma$ with $\hbox{CRdim } \Sigma=p-1$ and $2p-2 \leq \hbox{dim }
\Sigma \leq 2p+q-2$. Hence our Theorem 2 is far from being a corollary
of Theorem 3, even in the hypersurface case.

Instead of considering continuous CR functions or CR distributions,
one is led in the study of removable singularities in complex analysis
to consider functions holomorphic into wedges without growth condition
at all. Thus another notion of removability is as follows. By a wedge
{\it attached} to a generic manifold $M$, we mean an open {\it
connected} set which contains a wedge of edge $M$ at every point of
$M$ with a continuously varying direction in the normal bundle to
$M$. Let $\Phi\subset M$ be a proper closed subset of $M$. Call $\Phi$
{\it removable} if, given a wedge ${\cal W}_0$ attached to
$M\backslash \Phi$, there exists a wedge ${\cal W}$ attached to $M$
such that holomorphic functions into ${\cal W}_0$ extend
holomorphically into ${\cal W}$. (For a precise definition, see
Proposition 5.5).

We point out that the merit of our ``deformation philosophy'' is to
show that the two notions of removability are rather one and the
same. Theorems 1, 2 and 3 hold for both.

\medskip
\noindent
{\sc Theorem 4.} {\it Let $M, N$ be as in theorems 1,2,3 respectively
and let ${\cal W}_0$ be a wedge attached to $M\backslash N$. Then
there exists a wedge ${\cal W}$ of edge $M$ at $z_0$ such that
holomorphic functions into ${\cal W}_0$ extend holomorphically into
${\cal W}$.}

\smallskip
{\it Remark 1.} We give a proof of Theorem 4 with an everywhere
minimal manifold $M$, but it is possible to prove removability of $N$
in the wedge sense without any minimality assumption on the base
manifold.

\smallskip
{\it Remark 2.} Let $M$, $N$, $\Phi$ be as in Theorem 3. We choose a
formulation avoiding the notion of orbits, but we shall obtain a more
global result, if $N$ is assumed to be a closed embedded submanifold:
{\it $\Phi$ is fully removable if $\Phi$ does not contain any CR orbit
of $N$.}

\smallskip
{\it Remark 3.} The condition that the tangent space to $N$ does not
contain the full complex tangent space to $M$ at $z_0$ means, roughly,
that the CR geometry is not entirely absorbed at that point. Such a
hypothesis seems unavoidable and is commonly made in the hypersurface
case \cite{CS} \cite{JO2}. However, in the real analytic category, we
can dispense ourselves of it in Theorem 1 (Corollary 4.3) and
presumably also in Theorem 2.

\smallskip
{\it Remark 4.} We shall also derive in Theorem 5.2 below removability
of a proper closed subset $\Phi$ in a connected $C^1$ manifold $N$ of
codimension two in $M$ with the condition on the tangent spaces, for
general $p\geq 1$, thus extending Theorem 1. When $p\geq 2$ and $N$ is
not a $(p-1)$-CR manifold, the set $N\backslash N^{CR}$ of generic
points is nonempty; these points are naturally shown to be removable,
since one can prove there that CR functions on $M\backslash N$ are
locally uniformly approximable on compact subsets of $M\backslash N$
by holomorphic polynomials (Proposition 5.B) and then apply the
propagation method of Tr\'epreau or the one of Tumanov (both can be
applied). To finish out the proof of Theorem 2, one then invokes
Theorem 5.2 with the set $\Phi = N^{CR}\neq N$ of $(p-1)$-CR points of
$N$.

\smallskip
{\it Remark 5.} The content of the sufficient condition in Theorem 2
can be explained as follows, at least locally. For every holomorphic
function $g$ near $z_0$ with $g(z_0)=0$ and $\partial g(z_0) \neq 0$,
the complex hypersurface $\Sigma_g=\{g=0\}$ intersects $M$ along a
two-codimensional submanifold $N=\Sigma_g \cap M$, whenever $\Sigma_g$
is transversal to $M$ in $\C^n$ at $z_0$, {\em i.e.} $T_{z_0}
\Sigma_g+T_{z_0}M=T_{z_0}\C^n$. Then the restriction of $1/g$ to
$M\backslash N$ does not extend holomorphically into any wedge of edge
$M$ at $z_0$, since $T_{z_0}\Sigma_g$ absorbs the whole of the normal
bundle to $M$ at $z_0$. In that case, $N$ is a CR manifold with
$\hbox{CRdim} \ N=p-1$, since it is generic in $\Sigma_g$.

One can prove a converse statement when $N$ is CR and minimal.

Recall that an analytic wedge ${\cal W}^{an}$ with edge a CR manifold
$N$ is a complex manifold with edge $N$, smooth up to $N$, with
dimension equal to the rank of the bundle $TN+JTN$.

\smallskip
\noindent
{\sc Theorem 5.} {\it Let $M$ be a $C^{2,\alpha}$-smooth generic
manifold in $\C^n (n\geq 3)$ with $\hbox{CRdim }M = p \geq 2$ and let
$z_0\in M$. Let $N\ni z_0$ be a $C^{2,\alpha}$ CR submanifold of $M$
with $\hbox{codim}_M N =2$, $\hbox{CRdim } N= p-1$,
$T_{z_0}N+T_{z_0}^cM=T_{z_0}M$ which is minimal everywhere. Then there
exists a wedge ${\cal W}$ of edge $M$ at $z_0$ and an analytic wedge
${\cal W}^{an}$ with edge $N$ which is a closed complex hypersurface
in ${\cal W}$ such that every CR function on $M\backslash N$ extends
holomorphically into ${\cal W}\backslash {\cal W}^{an}$.}

\smallskip
Here is a description of the content of the paper. In Sections 1 and
2, we recall the notion of the defect of an analytic disc attached to
a generic manifold. This will be useful to insure the existence of a
good disc. In Section 3, we introduce isotopies of analytic discs and
delineate a continuity principle. In Section 4 and 5, we complete the
proof of Theorems 1 and 2 by using normal deformations of analytic
discs in the main Proposition 4.1.

\smallskip
{\it Acknoweledgment.} The author wishes to thank J. M. Tr\'epreau for
having helped him to considerably improve the results, their statement
and the overall quality of the article.

\smallskip
\smallskip
\noindent
{\bf 1. Preliminaries.} This paragraph is concerned with results of
Tumanov \cite{TU1} and extracted from Baouendi, Rotshchild and
Tr\'epreau \cite{BRT}.

Let $M$ be a CR generic submanifold of the complex euclidean space,
$T_{z_0}M$ its tangent space at $z_0 \in M$ and $T_{z_0}^cM=T_{z_0}M
\cap JT_{z_0}M$ its complex tangent space, where $J$ denotes the
standard complex structure on $T\C^n$. Set $p=\hbox{CRdim} \ M$ and
$q=\hbox{codim} \ M$ and assume that $p>0$ and $q>0$. In the
following, $M$ will always be of smoothness class $C^{2, \alpha}$. An
{\it analytic disc} in $\C^n$ is a continuous mapping
$A:\overline{\Delta}\to \C^n$ which is holomorphic in $\Delta$, where
$\Delta$ is the open unit disc in $\C$, $\Delta=\{\zeta\in \C; \
|\zeta|<1\}$, $\overline{\Delta}=\Delta\cup b\Delta$ and $b\Delta$ is
the unit circle in $\C$. We say that $A$ is {\it attached to $M$
through $z_0$} if $A(b\Delta)\subset M$ and $A(1)=z_0$. We shall
always assume all analytic discs to be in the Banach space ${\cal
B}=(C^{1,\alpha}(\overline{\Delta})\cap {\cal O}(\Delta))^n$ and we
set ${\cal B}_0=\{A\in {\cal B}; \ A(1)=z_0\}$ or ${\cal
B}_0=(C^{1,\alpha}_0 (\overline{\Delta})\cap {\cal O}(\Delta))^n$,
where $C^{1,\alpha}_0(\overline{\Delta})$ denotes the set of
$C^{1,\alpha}$ functions on $\overline{\Delta}$ vanishing at $1$.

Note with $z_0$ the constant disc and consider, for $\varepsilon>0$
the neighborhood \begin{equation} {\cal B}_{z_0, \varepsilon} = \{A
\in (C^{1,\alpha} (\overline{\Delta})\cap {\cal O} (\Delta))^n; \
A(1)=z_0, ||A-z_0||_{C^{1,\alpha}} <\varepsilon\} \subset {\cal B}
\end{equation} If $r \in C^{2, \alpha}(U, \R^q)$ is a definig function
for $M$ in a neighborhood $U$ of $z_0$ in $\C^n$, we introduce the
mapping $R: {\cal B}_{z_0, \varepsilon} \to {\cal F}$ defined by
\begin{equation} A \mapsto (\zeta \mapsto r(A(\zeta))), \end{equation}
where ${\cal F}=C^{1,\alpha}_0(b\Delta, \R^q)$ is the subspace
consisting of $\R^q$-valued functions of class $C^{1,\alpha}$ on $b
\Delta$ vanishing at $1$.

With these notations, the subset ${\cal A} \subset {\cal B}_{z_0,
\varepsilon}$ of discs that are attached to $M$ is given by
\begin{equation} {\cal A} = {\cal A}_{z_0,\varepsilon}= \{A \in {\cal
B}_{z_0, \varepsilon}; \ R(A) =0 \}. \end{equation}

\noindent
{\sc Lemma 1.1.} {\it ${\cal A}$ is a Banach submanifold of ${\cal
B}_{z_0, \varepsilon}$ parameterized by $(C^{1,
\alpha}_0(\overline{\Delta}) \cap {\cal O}(\Delta))^p$.}

\smallskip
{\it Proof.} By virtue of the implicit function theorem in Banach
spaces, it is sufficient to prove that the differential mapping
$R'(z_0): {\cal B}_0 \to {\cal F}$ possesses a continuous right
inverse $S$.

The matrix $r_z=(\frac{\partial r_j}{\partial z_k}(z_0))_{1\leq j\leq
q, 1\leq k \leq n}$ has rank $q$, since $M$ is generic. Let $D$ be a
$n\times q$ matrix such that $r_z(z_0)D=I_{q\times q}$. For
$\widehat{A} \in {\cal B}_0$ and $A \in {\cal B}_{z_0, \varepsilon}$
we have \begin{equation} [R'(A)\widehat{A}](\zeta)= 2 \hbox{Re} \
(r_z(A(\zeta))\widehat{A}(\zeta)), \end{equation} and if we set, for
$f \in {\cal F}$ \begin{equation} S(f)=\frac{1}{2} D(f+iT_1 f),
\end{equation} where $T_1$ denotes the Hilbert transform of $f$
vanishing at $1$, we get

\smallskip
 \hspace{2.5cm} $2R'(z_0)S(f) =
 2r_z(z_0)S(f)+2\overline{r_z(z_0)S(f)}$

 \hspace{2.5cm} $\ \ \ \ \ \ \ \ \ \ \ \ \ \ \ \ = r_z(z_0)D(f+iT_1f)+
 r_{\bar{z}}(z_0)\overline{D}(f-iT_1f)$

 \hspace{2.5cm} $\ \ \ \ \ \ \ \ \ \ \ \ \ \ \ \ =2f.$

\smallskip
\noindent
$S: {\cal F} \to {\cal B}_0$ is the continuous right inverse of
$R'(z_0)$ we searched for. Furthermore, the closed subspace $Ker \
R'(z_0)=E_0$ which parameterizes ${\cal A}$ in a neighborhood of $z_0$
consists of discs $\widehat{A} \in {\cal B}_0$ such that $2 \hbox{Re }
r_z(z_0) \widehat{A} \equiv 0$ on $b \Delta$. $E_0$ is isomorphic to
$(C^{1, \alpha}_0(\overline{\Delta}) \cap {\cal O}(\Delta))^p$, since
$rg \ r_z(z_0)=q$ and, for every function $f$, holomorphic in $\Delta$
and continuous up to $b\Delta$ such that $2 \hbox{Re } f \equiv 0$ on
$b \Delta$, there exists a constant $c \in \R$ such that $f\equiv ic$
in $\overline{\Delta}$ (if $f(1)=0$, necessarily $c=0$). This
completes the proof.

\smallskip
\noindent
{\sc Defect of an analytic disc.} Let $M$ be generic as before. We
shall define the defect of an analytic disc $A\in {\cal A}$ as
follows. Let $\Sigma(M) \subset \Lambda^{1,0}(T^{*}\C^n)$ denote the
conormal bundle to $M$ in $\C^n$. Its fiber at $z\in M$ are given by
$$\Sigma_z(M)=\{\omega \in \Lambda^{1,0}_z(T^{*}\C^n); \ \hbox{Im }
<\omega, X>= 0 \ \forall X \in T_z M \}.$$

We consider analytic discs $B : \overline{\Delta} \to \Lambda^{1,0}
\C^n$ attached to $\Sigma(M)$. Choosing
$(z_1,...,z_n,\mu_1,...,\mu_n)$ as holomorphic coordinates on
$\Lambda^{1,0}\C^n$, $B(\zeta)=(A(\zeta), \mu(\zeta))$ is attached to
$\Sigma(M)$ if and only if \begin{equation} A(b\Delta)\subset M \ \ \
\ \ \ \ \ \hbox{and} \ \ \ \ \sum_j \mu_j(\zeta)dz_j \in
\Sigma_{A(\zeta)} (M) \ \ \ \ \ \zeta \in b \Delta. \end{equation} The
set $V_A$ of discs that are attached to $\Sigma(M)$ can be equipped
with a vector space structure on the fiber component. For $\zeta \in
b\Delta$, we set \begin{equation} V_A(\zeta)= \{ \omega \in
\Sigma_{A(\zeta)}M; \ \omega= \sum_{j=1}^{n} \mu_j(\zeta)dz_j, \
(A(\zeta), \mu(\zeta)) \in V_A \}. \end{equation}

\smallskip
\noindent
{\sc Definition 1.2.} {\it If $\zeta \in b\Delta$ and $A$ is a disc
attached to $M$, the defect $def_{\zeta} A$ of $A$ at $\zeta$ is the
dimension of the vector space $V_A(\zeta)$.}

\smallskip
The canonical identification between $\Lambda^{1,0}\C^n$ and
$T^{*}\C^n$ enables one to identify also $\Sigma(M)$ to the
characteristic bundle of the CR vector fields on $M$, $(T^cM)^{\bot}
\subset T^{*} M$ (see {\cite{BRT} or \cite{ME}). Let ${\R^q}^{*}$
denote the dual space to $\R^q$.

\smallskip
\noindent
{\sc Proposition 1.3.} {\it If $\varepsilon >0$ is sufficiently small
and $A \in {\cal A}_{z_0, \varepsilon} = {\cal A}$,} the defect
$def_{\zeta} A$ is independent of $\zeta \in b \Delta$. {\it More
precisely,,
$$V_A(\zeta_0)=\{\xi \in T_{A(\zeta_0)}^{*}M, \ \xi = ib\nu
(\zeta_0)\partial r(A(\zeta_0)); \ \ \ \ \ \ \ \ \ \ \ \ \ \ \ \ \ \ \
\ \ \ $$
$$ \ \ \ \ \ \ \ \ \ \ \ \ \ \ \ \ \ \ \ \ \ b \in {\R^q}^{*} \hbox{
and } \zeta \mapsto b \nu (\zeta) r_z (A(\zeta)) \hbox{ extends
holomorphically to } \Delta\},$$ where $\nu(\zeta)$ is the unique
$q\times q$ invertible matrix with coefficients in $C^{1,
\alpha}(b\Delta, \R)$ such that $\nu(1)=I_{q\times q}$ and $\zeta
\mapsto \nu(\zeta)r_z(A(\zeta))D$ extends holomorphically to
$\Delta$.}

\smallskip
{\it Proof}. A disc $B(\zeta)=(A(\zeta),\mu(\zeta))$ is attached to
$\Sigma(M)$ if and only if \begin{equation} \sum_j \mu_j(\zeta)dz_j =
i t(\zeta)\partial r (A(\zeta))
\end{equation}
with $\mu_j(\zeta)$ extending holomorphically to $\Delta$ and $b
\Delta \ni \zeta \mapsto t(\zeta) \in {\R^q}^{*}$ being of class
$C^{1, \alpha}(b\Delta)$. As a consequence, the mapping
\begin{equation}
\zeta \mapsto i t(\zeta) r_z(A(\zeta)) D \in \C^q
\end{equation}
extends holomorphically to $\Delta$. The existence and uniqueness of a
matrix $\nu(\zeta)$ such that $\zeta \mapsto \nu(\zeta)r_z(A(\zeta))D$
extends holomorphically to $\Delta$ can be established, using
elementary Banach space techniques, by noting that $r_z(A(\zeta))D$ is
close to the identity if the size $\varepsilon$ of $A$ is sufficiently
small.

\noindent
Then $t(\zeta)r_z(A(\zeta))D (\nu(\zeta)r_z(A(\zeta))D)^{-1}
=t(\zeta)\nu(\zeta)^{-1}$ extends holomorphically to $\Delta$. Since
$t$ and $\nu$ are real and $\nu(1)=I_{q\times q}$,
$t(\zeta)=t(1)\nu(\zeta)$ for $\zeta \in b \Delta$ and the form for
$V_A(\zeta_0)$ follows if one sets $b=t(1)$.

The proof of Proposition 1.3 is complete.

\smallskip
For a fixed $\zeta_0 \in b \Delta$, we define the {\it evaluation map}
 ${\cal F}_{\zeta_0}: {\cal B} \to \C^n$ given by
$${\cal F}_{\zeta_0}(A)=A(\zeta_0).$$ Also, for $A \in {\cal B}$, the
{\it tangential direction} mapping at $\zeta=1$,
$$A \mapsto {\cal G}(A)= \frac{\partial A}{\partial \theta} (1).$$ If
${\cal F}_{\zeta_0}$ and ${\cal G}$ are restricted to the Banach
submanifold ${\cal A}$ constructed above, their differentials at $A
\in {\cal A}$ are linear applications
$${\cal F}_{\zeta_0}'(A): T_A{\cal A} \to T_{A(\zeta_0)} M \ \ \ \ \ \
\ \ \ {\cal G}'(A): T_A {\cal A} \to T_{z_0} M.$$ Notice that ${\cal
F}_1'(A)\equiv 0$, since ${\cal F}_1(A)= z_0$ for each $A \in {\cal
A}$.

\smallskip
\noindent
{\sc Theorem 1.4.} {\it Let $\zeta_0 \in b \Delta, \zeta_0 \neq 1$. If
 $\varepsilon$ is sufficiently small, $A \in {\cal A}={\cal A}_{z_0,
 \varepsilon}$ and $V_A(\zeta_0) \subset T_{A(\zeta_0)}^{*}M$ is
 defined by $(7)$,

 $(i) \ {\cal F}_{\zeta_0}'(A)T_A {\cal A} = V_A(\zeta_0)^{\bot} \ \ \
 \ \ \ \subset T_{A(\zeta_0)} M$

 $(ii) \ {\cal G}'(A)T_A {\cal A} = V_A(1)^{\bot} \ \ \ \ \ \ \ \
 \subset T_{z_0} M,$

\noindent
where orthogonality is taken in the sense of duality between $TM$ and
$T^{*}M$.}

\smallskip
Since $V_A(\zeta_0)$ can be viewed as a subspace of
$(T_{A(\zeta_0)}^cM)^{\bot} \subset T_{A(\zeta_0)}^{*}M$, we have the
following corollary of Theorem 1.4.

\noindent
{\sc Corollary 1.5.} {\it Under the assumptions of Theorem 1.4, the
codimension of ${\cal F}_{\zeta_0}'(A)T_A{\cal A}$ in
$T_{A(\zeta_0)}M$ and the codimension of ${\cal G}'(A)T_A{\cal A}$ in
$T_{z_0}M$ coincide and are equal to the defect of $A$. Furthermore,
the following inclusion holds
$$T_{A(\zeta_0)}^cM \subset {\cal F}_{\zeta_0}'(A)T_A{\cal A} \ \
\hbox{ and } \ \ T_{z_0}^cM \subset {\cal G}'(A)T_A{\cal A}.$$}

\vspace{-0.2cm}
\noindent
For the proof of Theorem 1, which is due in substance to Tumanov, see
Baouendi, Rothschild and Tr\'epreau \cite{BRT}. In Section 2 below, we
plain to use Corollary 1.5 and the independence of the defect of a
disc on everywhere minimal manifolds.

\smallskip
{\it Remark.} A better regularity of discs attached to
$C^{2,\alpha}$-smooth manifolds can be obtained by analysing Bishop's
equation as, for example, in \cite{TU2}.

\smallskip
\noindent
{\bf 2. Existence of a disc.} In sections 4 and 5, we shall need an
embedded disc $A \in {\cal A}$ such that $A(1)\in N$, $A(b \Delta
\backslash \{1\}) \subset M\backslash N$ and
$\frac{d}{d\theta}|_{\theta=0}A(e^{i\theta})\not\in T_{A(1)} N$ to
perform removal of singularities along that disc. In the present
section, we show how we can derive the existence of such a disc from
the hypothesis that $M$ is minimal at every point.

In the rest of the paper, the word neighborhood always means {\it
open} neighborhood.

Let $N\subset M$ be a $C^1$ submanifold through $z_0$ with
$\hbox{codim}_M N \geq 2$, $T_{z_0} N \not\supset T_{z_0}^cM$ and pick
a $C^1$ generic manifold $M_1 \subset M$ containing $N$ with $codim_M
M_1=1$. We can assume that $z_0=0$ is our reference point in a
coordinate system $(w,z)=(w_1,...,w_p,z_1,...,z_q), w=u+iv\in \C^p,
z=x+iy\in \C^q$ such that $T_0^cM=\C_w^p \times \{0\}$, $T_0 M= \C_w^p
\times\R^q_x$ and $T_0M_1=\{y=0, v_1=0\}$, $T_0^cM\cap T_0M_1=\{z=0,
v_1=0\}$. Then $M$ is given by the $q$ smooth scalar equations
\begin{equation} y=h(w,x) \ \ \ \ \ h \in C^{2,\alpha}, \ \ \ \ \
h(0)=dh(0)=0. \end{equation} Using the solution of Bishop's equation
given in \cite{TU2}, one can see that for sufficiently small $c>0$,
the embedded disc $A_c: \zeta \mapsto
A_c(\zeta)=(w_c(\zeta),z_c(\zeta))$ with $w$-component
$w_c(\zeta)=(c(1-\zeta),0,...,0)$ and $z$-component satisfying the
functional relation (called Bishop's equation)
$$x_c(\zeta)=-[T_1(h(w_c(.),x_c(.)))](\zeta) \ \ \ \ \ \zeta \in b
\Delta,$$ meets $M_1$ at two points exactly along its boundary, namely
$A_c(1)=0$ and another, say $A_c(e^{i\theta_c}),$ where $e^{i\theta_c}
\in b\Delta$ is close to $-1$. ($T_1$ denotes the Hilbert transform on
the unit circle normalized by the condition $(T_1u)(1)=0$,
$T_1u=Tu-u(1)$.)

Indeed, for every $\zeta\in b\Delta$, $\frac{d}{dc}|_{c=0} A_c(\zeta)
\in T_0^cM$, since $X_c(\zeta):=\frac{d}{dc}|_{c=0}x_c(\zeta)$
satisfies the functional relation $X=-T_1(2\hbox{Re}
((1-.)ch_{w_1}(0,X)))$ on $b\Delta$, hence $X\equiv 0$, by the
uniqueness in the solutions of Bishop's equation, and $Y=T_1X=0$
also. According to \cite{TU2}, given $\beta < \alpha$ with $\beta >0$,
there exists $c_0 >0$ such that for $c\leq c_0$, the mapping
$(c,\zeta)\mapsto A_c(\zeta)$ is $C^{2, \beta}$-smooth. Since
$\frac{d}{dc}|_{c=0} \frac{d}{d\theta}|_{\theta=0} z_c(e^{i\theta})=0$
and $\frac{d}{dc}|_{c=0} z_c(1)\equiv 0$, we obtain $|\frac{\partial^2
z_c}{\partial c\partial\theta} (\zeta)| \leq C(c^{\beta}
+|1-\zeta|^{\beta})$ for some constant $C>0$ depending on the second
derivatives of $h$ at $0$, hence $|\frac{dz_c}{d\theta} (e^{i\theta})|
\leq 3cC$ and also $|z_c(\zeta)|\leq
C(c|1-\zeta|(c^{\beta}+|1-\zeta|^{\beta}))$. As a consequence, one can
realize $A_c(b\Delta)$ as a graph over $w_c(b\Delta),
A_c(b\Delta)=\{(w, \psi(w)); \ w\in w_c(b\Delta)\}$ where $\psi$ is
$C^{2, \beta}$-smooth and satisfies $|\frac{d\psi}{dw}(w)| \leq
3C$. Since $M_1$ has codimension one in $M$ and the boundary of the
disc $w_c$ in $\C^p_w$ meets $T_0^cM \cap T_0 M_1=\{z=0, v_1=0\}$ at
two points $w_c(1)$ and $w_c(-1)$, the lifting $bA_c$ of $bw_c$ to
$\C^n \cap M$ necessarily meets $M_1$ at exactly two points, $A_c(1)$
and a second one, $A_c(e^{i\theta_c})$ with $e^{i\theta_c}$ close to
$-1$ in $b\Delta$, if $c_0$ is sufficiently small in order that
$A_{c_0}(b\Delta)$ is contained in a neighborhood $U$ of $0$ in $M$
with $M_1 \cap U$ very close to $T_0M_1$.

{\it Remark.} The property that $bA_c$ meets $M_1$ at exactly two
points is stable under $C^{2, \beta}$ perturbations of $A_c$.

Since $N$ is contained in $M_1$, $bA_c$ can meet $N$ by at most two
points in $b\Delta$. If $bA_c$ meets $N$ only at $\zeta=1$, the disc
satisfies our requirement.

If not, we shall obtain a good disc by slightly perturbing $A_c$. For
$\delta>0$, let ${\cal V}(A_c, \delta)$ denote the ball of radius
$\delta$ and center $A_c$ in $(C^{1,\alpha}(\overline{\Delta}) \cap
{\cal O}(\Delta))^n$.

\smallskip
\noindent
{\sc Proposition 2.1.} {\it Assume that $M$ is minimal at every point
and $\hbox{codim}_M N\geq 2$. Then, for each $\delta>0$ and $\beta
<\alpha$, there exists a $C^{2,\beta}$-smooth embedded disc $A \in
{\cal A}\cap{\cal V}(A_c, \delta)$ such that $A(1)=z_0$,
$A(b\Delta\backslash \{1\}) \subset M\backslash N$,
$\frac{d}{d\theta}|_{\theta=0} A(e^{i\theta})=v_0\not\in T_{z_0}^cM$
and $v_0\not\in T_{z_0}N$.}

\smallskip
{\it Proof.} Choose a disc $A \in {\cal V}(A_c,\delta)$ of minimal
possible defect, $A\in {\cal A}$ attached to $M$. If $\delta$ is
sufficiently small, $bA$ will meet $M_1$ at only two points, $A(1)=0
\in M_1$ and some $A(\zeta_1)=p_1 \in M_1$, $\zeta_1$ being close to
$-1$ in $b\Delta$. Indeed, discs in ${\cal V}(A_c, \delta)\cap {\cal
A}$ are close to $A_c$ in $C^1$ norm. If $p_1\in M_1\backslash N$, we
are done. So, suppose that $p_1 \in N$. Since $M$ is minimal at the
point $p_1$ and $A$ is of minimal possible defect, $A$ is of defect
$0$. Indeed, according to Theorem 1.4, the evaluation map ${\cal
F}_{\zeta_1}$ is of constant rank $2p+q-def A$ in a neighborhood
${\cal U}$ of $A$ in ${\cal A}$. Then the image of ${\cal U}$ by
${\cal F}_{\zeta_1}$ is a submanifold $\Sigma$ of $M$ containing $p_1$
with $dim \ \Sigma=2p+q-def A$ for which $T_z\Sigma \supset T_z^cM$
for all $z \in \Sigma$ by Corollary 1.5. Since $M$ is minimal at
$p_1$, necessarily $def A=0$.

Shrinking ${\cal U}$, we can assume that all boundaries of discs in
${\cal U}$ also meet $M_1$ at exactly two points where their tangent
direction is not tangent to $M_1$. According to the above, there
exists a $(2p+q-1)$-dimensional manifold $T \subset {\cal U}$ through
$A$ such that ${\cal F}_{\zeta_1}$ is a diffeomorphism of $T$ onto a
neighborhood of $p_1$ in $M_1$. Since $N$ is a proper submanifold of
$M_1$, there exist many discs $A_t$ in $T\subset {\cal A}$ with
boundary satisfying $A_t(b\Delta\backslash \{1\}) \subset M\backslash
N$.

According to Theorem 1.4, one can also find such a disc $A_t$ with
$\frac{d}{d\theta}|_{\theta=0} A(e^{i\theta})\not\in T_{z_0}^c M$.

Finally, the $C^{2,\beta}$-smoothness of the mapping $(t,\zeta)
\mapsto A_t(\zeta)$ is a consequence of the estimates in \cite{TU2},
if we consider $(C^{2,\alpha}(\overline{\Delta}) \cap {\cal
O}(\Delta))^p$ instead of $(C^{1,\alpha}(\overline{\Delta})\cap {\cal
O}(\Delta))^p$ as a parameter space to ${\cal A}$. The proof of
Proposition 2.1 is complete.

\smallskip
\noindent
{\bf 3. Isotopies of analytic discs and a continuity principle.} Let
$T_M\C^n=T\C^n|_M /TM$ denote the normal bundle to $M$ and let
$\eta_0\in T_M\C^n[z_0]$. By a {\it wedge of edge $M$ at} $(z_0,
\eta_0)$, we mean an open set of the form
$${\cal W}={\cal W}(U,z_0)=\{z+\eta; z\in U, \eta \in C\},$$ where $C$
is a conic neighborhood of some {\it nonzero} representative of
$\eta_0$ in $T_{z_0}\C^n$, identified with $\C^n$. The definition is
independent, in the germs, of the choice of the representative, in the
sense that each wedge contains another wedge for another choice. By
definition, such a wedge ${\cal W}$ contains a neighborhood of $z_0$
in $\C^n$ if $\eta_0=0$ in $T_M\C^n[z_0]$.

The proof of Theorem 1 and Theorem 2 starts as follows. Recall that we
assumed that $M$ is minimal at every point. Then, according to the
theorem of Tumanov, CR functions on $M\backslash N$ are wedge
extendible at every point of $M\backslash N$. The following definition
will be convenient for our purpose.

\smallskip
\noindent
{\sc Definition 3.1.} {\it An open {\rm connected} set ${\cal W}_0$
will be called a {\rm wedge attached to $M$} if there exists a {\rm
continuous} section $\eta: M \to T_M\C^n$ of the normal bundle to $M$
and ${\cal W}_0$ contains a wedge ${\cal W}_z$ of edge $M$ at
$(z,\eta(z))$ for every $z$ in $M$.}

\smallskip
Applying if necessary the edge of the wedge theorem (\cite{AY}) at
points where the direction of extendibility varies discontinuously,
the hypothesis on $M$ implies that there exists a wedge ${\cal W}_0$
attached to $(M\backslash N)$ to which CR functions on $M\backslash N$
have a holomorphic extension.

{\it Heuristics.} The aim, to achieve the proof of the two theorems,
is to show that such holomorphic functions extend holomorphically into
a wedge of edge $M$ at $z_0$, under the various hypotheses on $N$.
Notice that we cannot a priori know wether CR functions on
$M\backslash N$ are approximable by holomorphic polynomials on small
compact subsets of $M\backslash N$. This is why we are forced to use
instead the continuity principle by first deforming $M$ into ${\cal
W}_0$.

Fix a function $f\in CR(M\backslash N)$. It extends into some open
wedge ${\cal W}_0$ attached to $M\backslash N$.  Using a
$C^{3}$-smooth partition of unity on $(M\backslash N)\cap V$, for some
neighborhood $V$ of $z_0$ in $M$, we can smoothly deform $M$ into
${\cal W}_0$ leaving $N$ fixed and then replace $f$ by the restriction
of its extension to the $C^{2,\alpha}$ deformation $M^d$ of $M$.
Indeed, smooth deformations $d$ of $M$ into ${\cal W}_0$ arbitrarily
close to $M$ in $C^{2,\alpha}$ norm and fixing $N$ are possible, since
${\cal W}_0$ has a continuously varying direction over $M\backslash
N$.  Then, instead of a function $f\in CR(M\backslash N)$, we get a
function $f$, holomorphic into a neighborhood $\omega$ ($\equiv {\cal
W}_0$) of $M^d\backslash N$ in $\C^n$. The aim will subsequently be to
prove that such holomorphic functions extend into a wedge of edge $M$
at $z_0$.

A natural tool in describing envelopes of holomorphy of general open
sets in $\C^n$ is the {\it continuity principle}. Our version is the
following.

If $E\subset \C^n$ is any set, let ${\cal V}(E,r)$ denote $\{z\in
\C^n; \ \hbox{dist}(E,z)<r\}$ where dist denotes the polycircular
distance. Then $B(z,r)={\cal V}(\{z\},r)$ is the polydisc of center
$z$ and radius $r$.

Note with ${\cal O}(\omega)$ the ring of holomorphic functions into
the open set $\omega$.

\smallskip
\noindent
{\sc Lemma 3.2.} {\it Let $\omega \subset \C^n$ be an open connected
set and let $A: \overline{\Delta} \to \C^n$ be an analytic disc such
that $c|\zeta -\zeta'|<|A(\zeta)-A(\zeta')|<C|\zeta-\zeta'|$ for some
constants $0<c<C$ and every $\zeta, \zeta'\in \overline{\Delta}$. Set
$r=\hbox{dist}(A(b\Delta), b \omega))$ and $\sigma=\frac{rc}{2C}$.
Then for every holomorphic function $f\in{\cal O}(\omega)$, there
exists a function $F\in {\cal O}({\cal V}(A(\overline{\Delta}),
\sigma))$ with $F=f$ on ${\cal V}(A(b\Delta), \sigma).$}

\smallskip
{\it Proof.} For $z_0=A(\zeta_0), \zeta_0\in \overline{\Delta}$, let
$\sum_{k\in \N^n} f_k(z-z_0)^k$ denote the germ at $z_0$ of the
converging Taylor series defining $f$ there.  Take $\rho >0$ with
$\rho <r$.  (3.2 is trivial, if $r=0$).  Then the maximum principle
and Cauchy's inequalities yield
$$|f_k|=\frac{1}{k!}|D^k f(z_0)| \leq \frac{1}{k!} \sup_{z\in bA}
|D^kf(z)| \leq \frac{M_{\rho}(f)}{\rho^k},$$ where $M_{\rho}(f)=\sup
\{|f(z)|; \ z\in {\cal V}(A(b\Delta), \rho)\}<\infty$.  This proves
that the Taylor series of $f$ converges in $B(z_0, r)$ for every $z_0
\in A(\overline{\Delta}),$ defining an element $f_{z_0, r}\in {\cal
O}(B(z_0,r))$.  Now, if $B(z, \sigma)\cap B(z',\sigma)\neq \emptyset,
z=A(\zeta), z'=A(\zeta')$ and $\tau\in [\zeta, \zeta'],$ then
$|\tau-\zeta|\leq |\zeta'-\zeta| <\frac{2\sigma}{c}$ and
$$|A(\tau)-A(\zeta)|<C|\tau-\zeta|<\frac{2C\sigma}{c}=r.$$ Therefore
$A([\zeta,\zeta']) \subset B(z,r)\cap B(z',r)$ and the two holomorphic
funcions $f_{z,r}$ and $f_{z',r}$ coincide on the connected
intersection, hence define a function holomorphic into $B(z,r)\cup
B(z',r)$.  This proves that the $f_{z,\sigma}=f_{z,r}|_{B(z,\sigma)}$
stick together in a well-defined holomorphic function into ${\cal
V}(A(\overline{\Delta}), \sigma)$.

The proof of Lemma 3.2 is complete.

\smallskip
\noindent
{\sc Proposition 3.3.} {\it Let $\omega$ be an open connected set and
let $A_s:\overline{\Delta} \to \C^n, b\Delta \to \omega, 0\leq s \leq
1$ be a continuous family of analytic discs such that, for some
constants $0< c_s <C_s$, $c_s |\zeta-\zeta'| <
|A_s(\zeta)-A_s(\zeta')| < C_s |\zeta-\zeta'|$.  Assume
$A_1(\overline{\Delta}) \subset \subset \omega$ and set
$r_s=\hbox{dist} (A_s(b\Delta), b\omega), \sigma_s=\frac{r_s
c_s}{2C_s}$.  Then, for every holomorphic function $f\in {\cal
O}(\omega),$ there exist functions $F_s\in {\cal O}({\cal
V}(A_s(\overline{\Delta}), \sigma_s))$ such that $F_s\equiv f$ on
${\cal V}(A_s(b\Delta), \sigma_s)$.}

\smallskip
{\it Proof.} Let $I_0 \subset [0, 1]$ denote the connected set of $s_0
\in [0, 1]$ such that the statement holds true for every $s_0 \leq
s\leq 1$. $1 \in I_0$, according to Lemma 3.2. Let $s_1 <s_0$ be such
that $A_{s_1}(\overline{\Delta})\subset {\cal
V}(A_{s_0}(\overline{\Delta}), \sigma_{s_0})$.  Since $F_{s_0}\equiv
f$ into ${\cal V}(A_{s_0}(b\Delta), \sigma_{s_0})$ and
$\hbox{dist}(A_{s_1}(b\Delta), b\omega)=r_{s_1}$, the Taylor series of
$f$ at points $A_{s_1}(\zeta_1)=z_1, \zeta_1\in b\Delta$, converges in
each polydisc $B(z_1, r_{s_1})$.  As in the proof of Lemma 3.2, this
proves that there exists a function $F_{s_1}\in {\cal O}({\cal
V}(A_{s_1}(\overline{\Delta}), \sigma_{s_1}))$ with $F_{s_1} \equiv f$
in ${\cal V}(A_{s_1}(b\Delta), \sigma_{s_1})$.  Thus $I_0$ is both
open and closed, hence $I_0=[0, 1]$.

The proof of Proposition 3.3 is complete.

\smallskip
Let $\Phi$ be a proper closed subset of $M$.  Then $M\backslash \Phi$
is a smooth generic manifold, which will play, in the sequel, the role
of the $M$ in Definition 3.4 below.

\smallskip
\noindent
{\sc Definition 3.4.} {\it Let $M$ be generic. An embedded analytic
disc $A$ attached to $M$ is said to be} \ analytically isotopic to a
point in $M$ {\it if there exists a $C^1$-smooth mapping
$(s,\zeta)\mapsto A_s(\zeta), 0\leq s \leq 1, \zeta \in
\overline{\Delta}$, such that $A_0=A$, each $A_s$ is an embedded
analytic disc attached to $M$ for $0\leq s <1$ and $A_1$ is a constant
mapping $\overline{\Delta} \to \{pt\} \in M$.}

\smallskip
In the next section, to derive a proof of Theorem 1, we shall use
isotopies of embedded analytic discs without needing a particular
control of the size of the open neighborhood arising in our continuity
principle Proposition 3.3.

\smallskip
\noindent
{\sc Corollary 3.5.} {\it Let $M$ be generic, $C^{2,\alpha}$, let
$\Phi$ be a proper closed subset of $M$ and let $\omega$ be a
neighborhood of $M\backslash \Phi$ in $\C^n$. If an embedded disc $A$
attached to $M\backslash \Phi$ is analytically isotopic to a point in
$M\backslash \Phi$, there exists a connected open neighborhood
$\omega_A$ of $A(\overline{\Delta})$ in $\C^n$ such that for every
holomorphic function $f\in {\cal O}(\omega)$ there exists a
holomorphic function $f_A\in {\cal O}(\omega_A)$ with $f_A\equiv f$ in
a neighborhood of $A(b\Delta)$ in $\C^n$.}

\smallskip
Finally, the following proposition gives an auxiliary tool in place of
the approximation theorem of Baouendi and Treves.

\smallskip
\noindent
{\sc Proposition 3.6.} {\it Let $M, \Phi$ be as above.  Assume that
$M\backslash \Phi$ has the wedge extension property and let $z_0\in b
\Phi$.  Let $A\in {\cal A}_{z_1,\varepsilon}$ for some $z_1\in
M\backslash \Phi$, $|z_1-z_0|<<\varepsilon$ be an embedded disc
attached to $M\backslash \Phi$.  If $A$ is analytically isotopic to a
point in $M\backslash \Phi$, each $A_s$ of the isotopy being of size
$\varepsilon$, then $f\circ A|_{b\Delta}$ extends holomorphically to
$\Delta$ for every function $f\in CR(M\backslash \Phi)$.}

\smallskip
{\it Proof.}  Fix a function $f\in CR(M\backslash \Phi)$.  It extends
holomorphically and continuously into some wedge ${\cal W}_0$ attached
to $M\backslash \Phi$.  Since the isotopy is small, we can follow it
on small deformations $M^d$ of $M$ into ${\cal W}_0$ that fix
$\Phi$. Indeed, the solvability of the equation $R(A)=0$ as in Lemma
1.1 is stable under perturbations of $R$, since $R$ is a submersion.
Moreover, we may assume that $M^d$ depends on a small real parameter
$d\geq 0$, $M^0=M$ and the solution $(A_s)^d$ depends smoothly on
$(s,d)$.  Indeed, an implicit function theorem with parameters is
valid on Banach spaces.  If the deformation $d$ is sufficiently small,
this yields an analytic isotopy $(A_s)^d$ of $A^d$ to a point in $M^d
\backslash \Phi$.  According to Proposition 3.3, $f\circ
A^d|_{b\Delta}$ then extends holomorphically to $\Delta$. Letting $d$
tend to $0$, we get the desired result by the continuity of $f$ on
${\cal W}_0\cup (M\backslash \Phi)$ and the smoothness of $d\mapsto
A^d$.

The proof of Proposition 3.6 is complete.

\smallskip
Let $A_c$ denote the disc introduced in Section 2, $c$ small. In
Propositions 3.7 and 3.8 below, we give conditions which insure that
discs close to $A_c$ and attached to $M$ minus a singularity set $N$
are analytically isotopic to a point in $M\backslash N$.

\smallskip
\noindent
{\sc Proposition 3.7.} {\it Let $M$ be generic, $C^{2,\alpha}$-smooth,
let $N$ be a $C^1$-smooth proper submanifold of $M$ with
$\hbox{codim}_M N =3$, let $z_0\in N$ and assume that
$T_{z_0}N\not\supset T_{z_0}^cM$.  If $c>0$ is small enough, if $A_1
\in {\cal V}(A_c, \delta)\cap {\cal A}$ for $\delta>0, \delta<<c$
small enough satisfies $A_1(b\Delta\backslash \{1\}) \subset
M\backslash N$, there exists $\delta' >0$, $\delta' <<\delta$ such
that each disc $A\in {\cal V}(A_1,\delta')$ attached to $M\backslash
N$ is analytically isotopic to a point in $M\backslash N$.}

\smallskip
More generally, the following proposition holds.

\smallskip
\noindent
{\sc Proposition 3.8.} {\it Let $M$ be generic, $C^{2,\alpha}$-smooth,
let $N$ be a $C^1$-smooth proper submanifold of $M$ with
$\hbox{codim}_M N =2$, let $z_0\in N$ and assume that
$T_{z_0}N\not\supset T_{z_0}^cM$.  Let $\Phi$ be a proper closed
subset of $N$ with $z_0 \in b \Phi$.  If $c>0$ is small enough, if
$A_1 \in {\cal V}(A_c, \delta)\cap {\cal A}$ for $\delta>0, \delta<<c$
small enough satisfies $A_1(b\Delta\backslash \{1\}) \subset
M\backslash N$, there exists $\delta' >0$, $\delta' <<\delta$ such
that each disc $A\in {\cal V}(A_1,\delta')$ attached to $M\backslash
\Phi$ is analytically isotopic to a point in $M\backslash \Phi$.}

\smallskip
{\it Proof.} Suppose first that $T_0 N \cap \C_{w_1}= \{0\}$. We shall
show that the following stronger statement holds: if $c>0$ is
sufficiently small, there exists $\delta>0$, $\delta<<c$ such that
each disc $A\in {\cal V}(A_c,\delta)$ attached to $M\backslash \Phi$
is analytically isotopic to a point in $M\backslash \Phi$.

Let $\Gamma\subset M$ be a convex open conic neighborhood of the
positive $u_1$-axis with vertex $0$ such that $\overline{\Gamma}\cap
N=\{0\}$.  Choose a similar cone $\Gamma'$ with
$\overline{\Gamma'}\cap S^{2n-1}\subset\subset \overline{\Gamma}\cap
S^{2n-1}$, where $S^{2n-1}$ denotes the unit sphere in $\R^{2n}$.

If $c$ is sufficiently small, $A_c(b\Delta)\cap M_1$ consists of two
points, $0$ and a second one, say $p_1$ with $|p_1-p_1^0|<O(c^2)$,
where $p_1^0$ denotes the point in the $u_1$-axis with $u_1$
coordinate equal to $2c$ (in coordinates $(w,x)$ on $M$).  Take
$\delta >0, \delta<<c$ such that for each disc $A\in {\cal
V}(A_c,2\delta)\cap {\cal A}$, $A(b\Delta)\cap M_1$ also consists of
two points, $0$ and a second one in $\Gamma'$ close to $p_1^0$.  If
$\delta$ is small enough, we can furthermore insure that all discs in
${\cal V}(A_c, 2\delta)$ attached to $M$ are embedded discs and
intersect $M_1$ along their boundaries at two points, one of which is
$(2\delta)$-close to $0$ and a second one in $\Gamma$.

Fix a disc $A\in {\cal V}(A_c, \delta)$ attached to $M\backslash \Phi$
and assume that $A(1)\in M_1\backslash \Phi$ is $\delta$-close to $0$,
which can be done modulo an automorphism of $\Delta$. Write
$A(\zeta)=(w(\zeta), z(\zeta))$.  Since $N$ has codimension one in
$M_1$ and $\Phi$ does not disconnect $M_1$ near $0$, there exists a
$C^1$-smooth curve $\mu: [0,1] \to M_1\backslash \Phi$ with
$\mu(0)=A(1)\in M_1\backslash \Phi, \mu(1)\in \Gamma',
|\mu(s)|<\delta$ and $\mu(s)\in M_1\backslash \Phi$ for each $0\leq s
\leq 1$. Consider the disc $A_s(\zeta)=(w(\zeta)+\mu_w(s),
z_s(\zeta))$ where
$$x_s=-T_1h(w+\mu_w(s), x_s)+\mu_x(s) \ \ \ \ \ \hbox{on} \ \ b
\Delta.$$ Then $A_s(1)=\mu_s(1)$ and each $A_s$ is an embedded disc
belonging to ${\cal V}(A_c, 2\delta)$, hence attached to $M\backslash
\Phi$.  Therefore $A$ is analytically isotopic in $M\backslash \Phi$
with a disc $A_1=B$ such that $(w_1, x_1+ih(w_1,x_1))=B(1)\in \Gamma'
\cap M_1$ and $B(b\Delta)$ meets $M_1$ at another point in
$\Gamma$. Since the pure holomorphic $w$-component $w(\zeta)$ of $B$
is $(2\delta)$-close to $w_c(\zeta)=c(1-\zeta)$ in
$C^{1,\alpha}$-norm, and since $z_s$ is close to $z_c$ which satisfies
the estimates
$|z_c(\zeta)|=O(c|1-\zeta|(c^{\beta}+|1-\zeta|^{\beta}))$ and
$|\frac{dz_c}{d\theta} (e^{i\theta})=O(c)|$, the analytic isotopy
$B_s(\zeta)=((1-s)w(\zeta)+w_1, z_{B,s}(\zeta))$, where
$$x_{B,s}=-T_1h((1-s)w+w_1, x_{B, s})+x_1$$ for $0\leq s \leq 1$ joins
$B$ with the point $B(1)$ in $M\backslash \Phi$, each $B_s$ being an
embedded disc attached to $M\backslash N$, $0\leq s < 1$.  Indeed,
each $B_s$ also meets $M_1$ along its boundary at two points contained
in $ \Gamma$ for $0\leq s < 1$ according to the position of $M_1$ and
since, $w_1$ being close to $w_{c,1}$, $\{w_1(\zeta); \ \zeta \in b
\Delta\}\subset \C$ is contained in a set of the form $\{\hbox{Re}\
w_1 > -a |\hbox{Im} \ w_1|\}$, $a>0$ small.

Assume now that $T_0N \cap \C_{w_1}=\R_{u_1}$ and let $A_1\in {\cal
V}(A_c, \delta)\cap {\cal A}$ with $A_1(b\Delta\backslash \{1\})
\subset M\backslash N$, $A_1(1)=z_0$.  If $c$ and $\delta <<c$ are
sufficiently small, each disc $A\in {\cal V}(A_c, 2\delta)$ attached
to $M$ is an embedded disc whose boundary meets $M_1$ at two points
exactly. Since $N$ is one-codimensional in $M_1$, $N$ divides $M_1$ in
two connected components, $M^{+}_1$ and $M_1^{-}$. Therefore, we can
assume that $bA_1 \cap M_1$ consists of $A_1(1)=0$ and a point
$A_1(\zeta_1)=p_1\in M_1^{+}$.  Let $\delta'>0, \delta'<<\delta$ and
fix a disc $A\in {\cal V}(A_1, \delta')$ attached to $M\backslash
\Phi$, let $A(1)\in M_1\backslash \Phi$ be the point in $M_1\backslash
\Phi$ $\delta'$-close to $0$. Since $\Phi$ does not disconnect $M_1$
near $z_0$, there exists a $C^1$-smooth curve $\mu: [0, 1] \to M_1
\backslash \Phi$ with $\mu(0)=A(1), \mu(1)\in M_1^{+}$,
$|\mu(s)|<\delta'$ and $\mu(s)\in M_1\backslash \Phi$ for $0\leq s
\leq 1$. If $\delta'<<\hbox{dist} \ (p_1, N)$, perturbing the base
point of $A$ along the curve $\mu$ enables one to make an analytic
isotopy in $M\backslash \Phi$ of $A$ with an embedded disc $B$
attached to $M\backslash \Phi$ with the property that $bB \cap M_1$
consists of two points in $M_1^{+}$.

Then $B(\zeta)=(w(\zeta)+w_1, z(\zeta))$ with $w(\zeta)$ being
$(2\delta)$-close to $w_c(\zeta)$ in $C^{1, \alpha}$-norm, $w(1)=0$
and $w(\zeta)$ is the holomorphic $w$-component of $A$. Now we can
push away $B$ from $N$ by using a isotopy along a curve $\gamma:[0,
1]\to M_1^{+}$ such that $\gamma(s)$ belongs to a one dimensional
manifold $\Lambda \subset M_1$ with $\Lambda \cap N =\{z_1\}$ and
$T_{z_1} N + T_{z_1}\Lambda = T_{z_1} M_1$, some $z_1$ close to $B(1)$
in $M_1$.  $B_s(\zeta)$ is obtained as the perturbation of $B$ keeping
the same pure holomorphic $w$-component $w(\zeta)$ and satisfying
$B_s(1)=\gamma(s)$, $\gamma(0)=B(1)$.  In other words,
$B_s(\zeta)=(w(\zeta)+\gamma_w(s), z_s(\zeta))$ where
$$x_s=-T_1h(w+\gamma_w(s), x_s)+\gamma_x(s).$$ Differentiating the
equation with respect to $s$, it is possible to check that
$\frac{dB_s}{ds}(\zeta)$ is close to $\frac{d\gamma}{ds}(s)$ uniformly
in $\zeta$ for small $c, \delta, \delta'$ and then independent of
$\frac{dB_s}{d\theta}(e^{i\theta})$, which is close to be contained in
the $\C_{w_1}$-axis. Therefore the embedded disc $B_s$ goes away from
$N$ in the direction of $\Lambda\cap M_1^{+}$. Since $w$-embedded
discs which are sufficiently far away from $N$ are analytically
isotopic to a point in $M\backslash N$, $B$ and then $A$ are
analytically isotopic to a point in $M\backslash \Phi$.

The proof of Proposition 3.8 is complete.

\smallskip
\noindent
{\bf 4. Proof of theorem 1.}

The following proposition together with Proposition 2.1 and
Proposition 3.7 proves that the envelope of holomorphy of arbitrarily
thin neighborhoods of $M\backslash N$ in $\C^n$ always cover a wedge
of edge $M$ at every point of $N$, if $N$ is a three-codimensional
submanifold of $M$ as in Theorem 1.

\smallskip
\noindent
{\sc Proposition 4.1.} {\it Let $M$ be generic, $C^{2,\alpha}$-smooth,
let $z_0\in M$, let $N\ni z_0$ be a $C^1$ submanifold with $codim_M
N=3$ and $T_{z_0}N\not\supset T_{z_0}^cM$ and let $\omega$ be a
neighborhood of $M \backslash N$ in $\C^n$.  Assume there exists a
sufficiently small embedded analytic disc $A\in
C^{2,\beta}(\overline{\Delta})$ attached to $M$, $A(1)=z_0$, with $A(b
\Delta \backslash \{1\}) \subset M \backslash N$,
$\frac{d}{d\theta}|_{\theta=0} A(e^{i\theta})=v_0\not\in T_{z_0}^cM$,
$v_0 \not\in T_{z_0}N$ and all discs in ${\cal V}(A,\delta)$ attached
to $M\backslash N$ are analytically isotopic to a point in
$M\backslash N$, for some $\delta>0$ .  Then there exists a wedge
${\cal W}$ of edge $M$ at $(z_0, Jv_0)$ such that for every
holomorphic function $f \in {\cal O}(\omega)$ there exists a function
$F \in {\cal O}({\cal W})$ with $F=f$ in the intersection of ${\cal
W}$ with a neighborhood of $M\backslash N$ in $\C^n$.}

\smallskip
{\it Remark.} $M$ is not supposed to be minimal at any point.

\smallskip
{\it Proof.} Fix a function $f \in {\cal O}(\omega)$.  We shall
construct deformations of our given original disc $A$ as in \cite{TU3}
with boundaries in $M \cup \omega$ to show that the envelope of
holomorphy of $\omega$ contains a (very thin) wedge of edge $M$ at
$A(1)$. Instead of appealing to a Baouendi-Treves approximation
theorem, a version of which being not a priori valuable here, we shall
naturally deal with the help of the so-called continuity principle.

We can assume that $A(1)=0$ and that $M$ is given in a coordinate
system as in $(10)$.  Set $A(\zeta)=(w(\zeta), x(\zeta)+ i
y(\zeta))$. Since $T_0N \not\supset T_0^cM$ and $v_0=\frac{d}{d
\theta}|_{\theta=0} A(e^{i\theta}) \not\in T_0N$, we can choose a
$C^1$ generic one-codimensional submanifold $K \subset M $ in $\C^n$
containing $N$ with tangent space $T_0 K $ not containing $v_0$.
Making a unitary transformation in the $w$-space, we get that $T_0 K
\cap T_0^cM$ is given by $\{z=0, v_1=0\}$ and then making a complex
linear transformation in $\C^n$ stabilizing $\{y=0\}$, we can assume
that $K$ has equation $v_1=k(u_1,w_2,...,w_p,x)$, where $k(0)=0,
dk(0)=0$. Since $v_0\not\in T_0K$, the projection on the $v_1$-axis of
$v_0$ is nonzero.

Let $\mu=\mu(w,x)$ be a $C^{\infty },$ $\R$-valued function with
support near the point $(w(-1), x(-1))$ that equals 1 there and let
$\kappa: \R^{q} \to \R^{q}$ be a $C^{\infty }$ function with
$\kappa(0)=0$ and $\kappa'(0)=Id$.  We can assume that the supports of
$\mu$ and $\kappa$ are sufficiently concentrated in order that every
manifold $M_t$ with equation \begin{equation}
y=H(w,x,t)=h(w,x)+\kappa(t)\mu(w,x) \end{equation} is contained in
$\omega$ and the deformation is localized in a neighborhood of $A(-1)$
in $\C^n$.  Let $\chi= \chi(\zeta)$ be a smooth function on the unit
circle supported in a small neighborhood of $\zeta=-1$ that will be
chosen later.  For every small $t$, we shall consider the natural
perturbation $A_t$ of the given disc $A$ with the same holomorphic
$w$-component as $A$ that is attached to the union of the $M_t$.  More
precisely, $A_t(\zeta)=(w(\zeta), x_t(\zeta)+iy_t(\zeta))$ is obtained
by solving the following equation on the unit circle \begin{equation}
x_t= -T_1H(w, x_t, t\chi).  \end{equation} According to \cite{TU2},
the solution of $(12)$ exists and depends in a $C^{2, \beta}$ fashion
with respect to all variables, no matter $0 < \beta'$ is $< \beta$,
provided everything is very small. Rename $\beta'$ as $\beta$.  We let
$\Pi$ denote the canonical bundle epimorphism $\Pi: T\C^n |_M \to
T\C^n |_M / TM$ and consider the $C^{1,\beta}$ mapping
\begin{equation} D: \ \ \ \ \ \R^{q} \ni t \ \longmapsto \ \Pi \left(
- \frac{\partial A_t}{\partial \zeta} (1) \right) \in T_0 \C^n / T_0 M
\simeq \R^{q}.  \end{equation}

\medskip
\noindent
{\sc Lemma 4.A.} {\it $\chi$ can be chosen in order that $ \ rk \
D'(0)= q$.}

\medskip
{\it Proof.} We, for completeness, include a proof of this result,
originally due to Tumanov in \cite{TU3}. Set $t= (t_1,...,t_{q})$.
Differentiating the equation $y_t( \zeta) = H(w(\zeta), x_t(\zeta),
t\chi(\zeta)), \zeta \in b \Delta$, with respect to $t_j, j=1,...,q,$
we obtain that the holomorphic disc $\frac{\partial}{\partial t_j}
|_{t=0} A_t(\zeta) = \dot{A}(\zeta) = (0,
\dot{X}(\zeta)+i\dot{Y}(\zeta))$ satisfies the following equation on
the unit circle \begin{equation} \dot{Y}= H_x\circ A \dot{X} + \chi
H_{t_j} \circ A.  \end{equation} We also introduce some notations. For
a $C^{1, \beta}$-smooth function $g(\zeta)$ on the unit circle with
$g(1)=0$, we write
$${\cal J}(g)= \frac{1}{\pi} \int_0^{2\pi}
\frac{g(e^{i\theta})}{|e^{i\theta} - 1|^2} d \theta,$$ where the
integral is understood in the sense of principal value.  Then, if $g
\in C^{1, \beta}(\overline{\Delta})$ is holomorphic in $\Delta$ and
vanishes at $1$, we have \begin{equation} {\cal J}(g) = - \frac{
\partial g}{\partial \zeta}(1) = i \frac{d}{d \theta} {\mid}_{\theta =
0} g(e^{i \theta}).  \end{equation} Notice also that for $C^{1,
\beta}$ real-valued functions $g, g'$ with $g(1)=g'(1)=0$, applying
$(15)$ to the holomorphic function $(g+iT_1g)(g'+iT_1g')$ vanishing to
second order at $1$, we obtain \begin{equation} {\cal J}(gg' - T_1g
T_1g')=0.  \end{equation}

Associate with $A$ and $H$ a $q \times q$ matrix-valued function
$G(\zeta)$ on the unit circle as a solution to the equation
\begin{equation} G= I + T_1 (GH_x \circ A).  \end{equation} The
definition of $G$ implies that $G(1)=I$ and $T_1 G = -GH_x \circ A +
H_x \circ A(1)= -GH_x \circ A$, since $A(1)=0$ and $h(0)=0, dh(0)=0$.
Using $(14)$ and $G$, we can write on the unit circle

\smallskip

\hspace{4cm} $G \chi H_{t_j} \circ A = G(\dot{Y} - H_x \circ A
\dot{X})$

\hspace{6cm} $=G\dot{Y} - (T_1G)(T_1\dot{Y})$

\hspace{6cm} $=\dot{Y}+ (G-I)\dot{Y} -T_1(G-I) T_1 \dot{Y}.$

\smallskip
\noindent
By virtue of $(16)$, \begin{equation} {\cal J}(G\chi H_{t_j} \circ A )
 = {\cal J}(\dot{Y}).  \end{equation} On the other hand, according to
 $(15)$ and the fact that $\dot{X}=-T_1\dot{Y}$, \begin{equation}
 {\cal J}(\dot{Y})+i{\cal J}(T_1\dot{Y})= i\frac{\partial
 \dot{Z}}{\partial \zeta} (1)= \frac{d}{d\theta}|_{\theta=0}
 (\dot{X}+i\dot{Y}).  \end{equation} Identifying the imaginary part of
 the two extreme terms and taking $(14)$ into account, we have
 \begin{equation} {\cal J}(T_1\dot{Y})=\frac{d}{d\theta}|_{\theta=0}
 \dot{Y} = \frac{d}{d\theta}|_{\theta=0} (H_x\circ A \dot{X}+\chi
 H_{t_j} \circ A)=0, \end{equation} if we choose $\chi$ in order that
 $\chi$ is equal to zero near $\zeta=1$ and since $\dot{X}(1)=0$,
 $dH(0)=0$.  $(18)$, $(19)$ and $(20)$ therefore yield
$$i \frac{\partial \dot{Z}}{\partial \zeta}(1) ={\cal J}(G\chi H_{t_j}
\circ A).$$ Natural coordinates on $T_0\C^n / T_0 M $ being given by
$y_1,...,y_{q}$, we obtain in these coordinates \begin{equation}
\frac{\partial D}{\partial t_j} (0) = \Pi \left(- \frac{\partial
\dot{Z}}{\partial \zeta}(1) \right) = {\cal J}(G \chi H_{t_j} \circ
A).  \end{equation} Furthermore, choose $\chi$ in order that ${\cal
J}(\chi)=1$ and the support of $\chi$ is concentrated near $\zeta= -1$
so that the vectors ${\cal J}(G \chi H_{t_j}\circ A)$ are close to the
vectors $G(-1)H_{t_j}\circ A(-1)$ and linearly independent, for
$j=1,...,q$ respectively.  This is possible, since $G$ is non singular
at every point on the unit circle and the $H_{t_j}\circ A(-1)$,
$j=1,...,q$ are linearly independent by the choice of $\kappa$.  This
completes the proof of lemma 4.A.

We deduce:

\medskip
\noindent
{\sc Corollary 4.B.} {\it There exists a neighborhood ${\cal T}$ of
$0$ in $\R^{q}$ such that $\{s D(t); \ s>0, t \in {\cal T} \}$
contains $\Gamma_0$ an open cone with vertex $0$ in $\R^{q}$ and $\Pi
( -\frac{\partial A}{\partial \zeta}(1)) \in \Gamma_0$.}

\medskip
We shall need another type of deformations of our disc.  Let $\tau$
 denote a supplementary real parameter, $\tau_0 >0$ and ${\cal V}$ a
 neighborhood of $0$ in $\C^{p-1}$. For $t\in {\cal T}$,
 $|\tau|<\tau_0$, $a\in {\cal V}$, we consider the disc
 \begin{equation} A_{t,\tau,a}(\zeta)= (e^{i\tau}w_1(\zeta),w_2(\zeta)
 +a_2(\zeta-1),...,w_p(\zeta)+a_p(\zeta-1), x_{t,\tau,a}( \zeta)+i
 y_{t,\tau,a}(\zeta)) \end{equation} where $x_{t, \tau,a}$ is the
 $C^{2, \beta}$-smooth in all variables solution of the following
 equation on $b \Delta$ \begin{equation} x_{t, \tau,a}=
 -T_1H(e^{i\tau}w_1(\zeta),w_2(\zeta)+
 a_2(\zeta-1),...,w_p(\zeta)+a_p(\zeta-1), x_{t, \tau,a}, t \chi).
 \end{equation} Introduce also the smooth mapping \begin{equation} E:
 \ \ \ \ \ {\cal T} \times I_{\tau_0}\times {\cal V} \ni (t, \tau, a)
 \ \longmapsto \ \frac{d}{d \theta} |_{\theta=0} A_{t,
 \tau,a}(e^{i\theta}) \in T_0 M.  \end{equation} According to Lemma
 4.A, we have that $rk \ E'(0) = 2p+q-1$ for a good choice of the
 function $\chi$, since we defined $w_{t, \tau,a}(\zeta)$ not
 depending on $t$ and the partial rank of $E$ with respect to
 $(\tau,a)$ is equal to $(2p-1)$.  So is the rank of the mapping
 derived from $E$ by taking normalized tangent vectors with respect to
 a hermitian structure on $\C^n$.

Therefore we obtained

\smallskip
\noindent
{\sc Lemma 4.C.} {\it $\chi$ can be chosen in order that the following
holds: there exist $\tau_0 >0$, ${\cal T}$ a neighborhood of $0$ in
$\R^q$ and ${\cal V}$ a neighborhood of $0$ in $\C^{p-1}$ such that
the set \begin{equation} \hat{\Gamma}_0= \{s
\frac{dA_{t,\tau,a}}{d\theta} (1); \ s <0 \hbox{ or } s>0, t \in {\cal
T}, \tau\in I_{\tau_0}, a\in {\cal V}\} \end{equation} is a
$(2p+q)$-dimensional open connected bicone with vertex $0$ in the
$(2p+q)$-dimensional space $T_0M$ containing
$v_0=\frac{d}{d\theta}|_{\theta=0} A(e^{i\theta})$.}

\smallskip
\noindent
What we call a bicone is a union of two linear cones with same vertex
and symmetric opposite directions.  As a consequence, we can read
Lemma 4.C inside the base manifold as follows.

\smallskip
\noindent
{\sc Lemma 4.D.} {\it Shrinking the open neighborhoods ${\cal T},
I_{\tau_0}, {\cal V}$ if necessary, there exists an open connected arc
$1 \in I_1 \subset b \Delta$ such that the set \begin{equation}
\hat{\Gamma}= \{A_{t, \tau,a}(\zeta); \ t \in {\cal T}, \tau\in
I_{\tau_0}, a\in {\cal V}, \zeta \in I_1\} \end{equation} is a
$(2p+q)$-dimensional closed connected (nonlinear) truncated bicone
contained in $M$.}

\smallskip
\smallskip
We now introduce a third deformation of $A$, consisting in varying the
base point.  We let $A_{t, \tau, a, p_0}$ denote the disc attached to
$M \cup \omega$ which is the perturbation of $A_{t, \tau,a}$ that
passes through the point $p_0 \in K$ with coordinates
$(u_1^0,w_2^0,...,w_p^0, x^0)$ on $K$, {\em i.e.}  $A_{t, \tau,a,
p_0}(1)=p_0$ and {\small \begin{equation}
A_{t,\tau,a,p_0}(\zeta)=(e^{i\tau}w_1(\zeta)
+u_1^0+iv_1^0,w_2(\zeta)+a_2(\zeta-1)
+w_2^0,...,w_p(\zeta)+a_p(\zeta-1)+w_p^0, x_{t,\tau,a,p_0}( \zeta)+i
y_{t,\tau,a,p_0}(\zeta)) \end{equation}} where $x_{t,\tau,a,p_0}$ is
the solution of Bishop's equation with parameters \begin{equation}
x_{t,\tau,a,p_0}=-T_1H (e^{i\tau}w_1+u_1^0+iv_1^0,w_2+a_2(.-1)
+w_2^0,...,w_p+a_p(.-1)+w_p^0,x_{t,\tau,a,p_0},t\chi)+x^0,
\end{equation} and where $v_1^0 =
k(u_1^0,w_2^0,...,w_p^0,x_1^0,...,x_q^0)$, so $A_{t,\tau,a,p_0}(1)=p_0
\in K$.

For convenience, we shall allow us to shrink a finite number of times
the open set ${\cal T}\times I_{\tau_0}\times {\cal V} \times {\cal
K}$ {\it without explicit mention} in the rest of the proof of
Proposition 4.1.

Notice that if we let $\tau =0$ and $a=0$, the set of points $A_{t,
0,0, p_0}(\zeta)$ for $\zeta$ in the interior of a neighborhood
$\Delta_1$ of $1$ in $\overline{\Delta}$ always cover a wedge of edge
$M$ at $0$ in the direction $Jv_0$ when $t$ varies in a neighborhood
${\cal T}$ of $0$ in $\R^{q}$ and $p_0$ in a neighborhood ${\cal K}$
of $0$ in $K$, according to Corollary 4.B and the fact that
$\frac{d}{d\theta}|_{\theta=0}A(e^{i\theta}) \not\in T_0K$.  In the
non singular case, when a local approximation theorem is valid for CR
functions, this leads to propagation of wedge extendibility along a
disc, as in \cite{TU3}.

Since, however, our $f \in {\cal O}(\omega)$ possibly has
singularities on $N$, a different argument is needed.  Notice that
every embedded disc $A_{t,\tau,a,p_0}$ with sufficiently small
parameters $(t, \tau, a, p_0)$ is close to $A$ and is analytically
isotopic to a point in $\omega$ for $p_0\in {\cal K}\backslash N$.
Indeed, $A_{t,\tau,a,p_0}$ is first analytically isotopic in $\omega$
to $A_{0,\tau,a,p_0}$ by construction and $A_{0,\tau,a,p_0}$ is
attached to $M\backslash N$.  Clearly, since
$\frac{d}{d\theta}|_{\theta=0} A(e^{i\theta})=v_0$ is not tangent to
$N$, and $A=A_{0,0,0,0}$ satisfies $A(b\Delta\backslash \{1\})\subset
M\backslash N$, each disc $A_{0,\tau,a,p_0}$ with $p_0$ not in $N$
cannot meet $N$ along its boundary, provided all the parameters are
small. On the contrary, discs in the family with base point $p_0\in
{\cal K}\cap N$ meet $N$, but only at $p_0$.  Then the fact that for
$p_0\in {\cal K}\backslash N$, $A_{0,\tau,a,p_0}$ is analytically
isotopic to a point in $M\backslash N$ is a consequence of our
hypothesis that all discs in ${\cal V}(A, \delta)$ attached to
$M\backslash N$ are analytically isotopic to a point in $M\backslash
N$, if $\delta>0$ is very small.

Since our discs with base point $p_0 \in {\cal K}\backslash N$ are
analytically isotopic in $\omega$ to a point, our given function $f
\in {\cal O}(\omega)$ can be analytically extended in a neighborhood
$\omega_A$ of $A_{t,\tau,a, p_0}(\overline{\Delta})$ in $\C^n$, if we
make use of the so-called continuity principle along the isotopy
(Proposition 3.3). Precisely, there exist a neighborhood $\omega_A$ of
$A(\overline{\Delta})$ in $\C^n$ and $f_A \in {\cal O}(\omega_A)$ with
$f_A\equiv f$ in a neighborhood of $A_{t, \tau, a , p_0}(b \Delta)$.
On the contrary, discs meeting $N$ necessarily lack a similar
extension property.

We obtained that for each $p_0\in {\cal K}\backslash N$, every CR
function $f\in \hbox{CR} (M\backslash N)$ has the property that
$f\circ A_{t,\tau,a,p_0}|_{b\Delta}$ extends holomorphically into
$\Delta$.

Let $v_0\in\Gamma$ be a $q$-dimensional proper linear bicone in the
$(2p+q)$-dimensional space $T_0 M$ and contained in
$\widehat{\Gamma}_0$ such that the projection $T_0\Gamma \to T_0M/
T_0^cM$ is surjective and $\overline{\Gamma}\cap T_0^cM=\{0\}$.  Let
${\cal P}$ denote the set of parameters
$${\cal P}=\{(t,\tau,a)\in {\cal T}\times I_{\tau_0}\times {\cal V}; \
\frac{d}{d\theta}A_{t,\tau,a}(1) \in \Gamma\}.$$ Then ${\cal P}$ is a
$C^1$-smooth $(q-1)$-dimensional submanifold of ${\cal T}\times
I_{\tau_0}\times {\cal V}$.  Shrinking ${\cal T}, I_{\tau_0}, {\cal
V}, {\cal K}$ and $\Delta_1$ if necessary, the set
$${\cal W}=\{A_{t,\tau,a,p_0}(\zeta); \ (t,\tau,a)\in {\cal P}, p_0\in
{\cal K}, \zeta \in \stackrel{\circ}{\Delta}_1\}$$ contains a wedge of
edge $M$ at $0$. Furthermore, the mapping
$${\cal P}\times {\cal K}\times \stackrel{\circ}{\Delta}_1 \ni
(t,\tau,a,p_0,\zeta) \mapsto A_{t,\tau,a,p_0}(\zeta) \in
\C^n\backslash M$$ becomes a smooth embedding. Since the mapping
remains injective on ${\cal K}\backslash N$, we can set unambiguously
$$F(z):= \frac{1}{2i\pi} \int_{b\Delta} \frac{f\circ
A_{t,\tau,a,p_0}(\eta)}{\eta - \zeta} d\eta$$ as a value at points
$z=A_{t,\tau,a,p_0}(\zeta)$ for an extension of $f|_{M\backslash N}$,
$p_0 \in {\cal K}\backslash N, (t,\tau,a)\in {\cal P}, \zeta \in
\stackrel{\circ}{\Delta}_1$.  Since $f$ extends holomorphically to the
interior of these discs, we get a continuous extension $F$ on each
$A_{t,\tau,a,p_0}(\Delta_1)$, $p_0\in {\cal K}\backslash N$.  Thus,
the extension $F$ of $f|_{M\backslash N}$ also becomes continuous on
$$({\cal W}\backslash N_{\cal P}) \cup (M\backslash N),$$ where
$$N_{\cal P}=\{A_{t,\tau,a,p_0}(\zeta); \ (t,\tau,a)\in {\cal P},
p_0\in {\cal K}\cap N, \zeta \in \stackrel{\circ}{\Delta}_1 \}.$$
Since $f|_{M\backslash N}$ extends analytically to a neighborhood of
$A_{t, \tau, a, p_0}(\overline{\Delta})$, $F$ is holomorphic into
${\cal W}\backslash N_{\cal P}$. Indeed, fix a point $(\tilde{t},
\tilde{\tau}, \tilde{a}, \tilde{p}_0) \in {\cal P} \times ({\cal K}
\backslash N)$ and let $\tilde{\cal P} \times \tilde{\cal K}$ be a
neighborhood of $(\tilde{t}, \tilde{\tau}, \tilde{a}, \tilde{p}_0)$ in
${\cal P} \times ({\cal K}\backslash N)$ such that for each $(t, \tau,
a, p_0)\in \tilde{\cal P}\times \tilde{\cal K}$,
$A_{t,\tau,a,p_0}(\overline{\Delta})$ is contained in some
neighborhood $\tilde{\omega}$ of $A_{\tilde{t}, \tilde{\tau},
\tilde{a}, \tilde{p}_0}(\overline{\Delta})$ in $\C^n$ such that there
exists a holomorphic function $\tilde{f} \in {\cal O}(\tilde{\omega})$
with $\tilde{f}$ equal to $f$ near $A_{\tilde{t}, \tilde{\tau},
\tilde{a}, \tilde{p}_0}(b\Delta)$.  Let $\tilde{\zeta}\in
\stackrel{\circ}{\Delta}_1$ and $\tilde{z}=A_{\tilde{t}, \tilde{\tau},
\tilde{a}, \tilde{p}_0}(\tilde{\zeta})$. To check that the previously
defined function $F$ is holomorphic in a neighborhood of $\tilde{z}$,
we note that for $z=A_{t, \tau, a, p_0}(\zeta)$, $(t, \tau, a, p_0)
\in \tilde{\cal P}\times \tilde{\cal K}$, $\zeta$ in some neighborhood
$\tilde{\Delta}_1$ of $\tilde{\zeta}$ in $\stackrel{\circ}{\Delta}_1$,
$\tilde{f}(z)$ is given by the Cauchy integral formula
$$\tilde{f}(z)=\frac{1}{2i\pi} \int_{b \Delta} \frac{\tilde{f} \circ
A_{t, \tau, a ,p_0}(\eta)}{\eta-\zeta} d\eta= \frac{1}{2i\pi} \int_{b
\Delta} \frac{f\circ A_{t, \tau, a ,p_0}(\eta)}{\eta-\zeta} d\eta=
F(z).$$ As a consequence, $\tilde{f}(z)=F(z)$ for $z$ in a small
neighborhood of $\tilde{z}$ in $\C^n$, since the mapping $(t, \tau, a,
p_0,\zeta) \mapsto A_{t, \tau, a, p_0}(\zeta)$ from $\tilde{\cal
P}\times \tilde{\cal K} \times \tilde{\Delta}_1$ to $\C^n$ has rank
$2n$ at $(\tilde{t}, \tilde{\tau}, \tilde{a}, \tilde{p}_0,
\tilde{\zeta})$.

This proves that $F$ is holomorphic into ${\cal W}\backslash N_{\cal
P}$.
 
By shrinking $\omega$ near $0$, which does not modify the possible
disc deformations, we can insure that
$$\omega \cap {\cal W}$$ is connected, since $\overline{\Gamma}\cap
T_0^cM=\{0\}$ and then also
$$\omega \cap ({\cal W}\backslash N_{\cal P}),$$ since $N_{\cal P}$ is
a closed two-codimensional submanifold of ${\cal W}$. Therefore $f\in
{\cal O}(\omega)$ and $F\in {\cal O}({\cal W}\backslash N_{\cal P})$
stick together in a single holomorphic function in $\omega \cup ({\cal
W} \backslash N_{\cal P})$, since both are continuous up to
$M\backslash N$, which is a uniqueness set, and coincide there.

However, the proof of Proposition 4.1 will not be finished until we
get rid of $N_{\cal P}$.  This is why we introduced the supplementary
parameters $\tau$ and $a$.

Since $\hat{\Gamma}_0$ is a $(2p+q)$-dimensional bicone in the
$(2p+q)$-dimensional space $T_0M$, we can choose two $q$-dimensional
bicones $v_0\in\Gamma_2$ and $\Gamma_2'$ contained in $\hat{\Gamma}_0$
with the properties that

\smallskip
\smallskip
\hspace{1.5cm}
\begin{minipage}[t]{12.5cm}
\noindent
$(i)$ The projections $T_0 \Gamma_2 \to T_0 M / T_0^cM$ and
$T_0\Gamma_2' \to T_0M /T_0^cM$ are surjective, $\overline{\Gamma_2},
\overline{\Gamma_2} ' \cap T_0^cM=\{0\}$.  Moreover, for $v_2 \in
\Gamma_2$ there exists $v_2' \in \Gamma_2'$ such that $v_2 - v_2' \in
T_0^cM/(T_0^cM\cap T_0N)\backslash \{0\}$ and vice versa.

\noindent
$(ii)$ \ There exist open bicones $\hat{\Gamma}_2$ and
$\hat{\Gamma}_2'$ in $\hat{\Gamma}_0$ with $\Gamma_2 \subset
\hat{\Gamma}_2, \Gamma_2' \subset \hat{\Gamma}_2'$ and $\hat{\Gamma}_2
\cap \hat{\Gamma}_2' = \emptyset.$
\end{minipage}

\smallskip
\smallskip
Set ${\cal P}_2= \{(t, \tau, a); \frac{d}{d \theta} A_{t, \tau, a}(1)
\in \Gamma_2 \}$ and similarly for ${\cal P}_2'$. Then condition $(i)$
insures that the two following wedges with edge $M$ at $(0, \eta_0)$,
$\eta_0=Jv_0 \ \hbox{mod} \ T_0M$, \begin{equation} {\cal W}_2=
\{A_{t, \tau, a, p_0}(\stackrel{\circ}{\Delta}_1); \ (t, \tau, a) \in
{\cal P}_2, p_0 \in {\cal K} \} \ \ \ \ \ {\cal W}_2' = \{ A_{t, \tau,
a, p_0}(\stackrel{\circ}{\Delta}_1); \ (t, \tau, a) \in {\cal P}_2',
p_0 \in {\cal K} \} \end{equation} contain a wedge ${\cal W}$ of edge
$M$ at $(0, \eta_0)$.  According to the above, one can construct two
holomorphic extensions $F_2$ and $F_2'$ of $f$ into ${\cal
W}_2\backslash N_{{\cal P}_2}$ and ${\cal W}_2'\backslash N_{{\cal
P}_2'}$ respectively, where $N_{{\cal P}_2}$ and $N_{{\cal P}_2'}$
denote the unattainable sets,
$$N_{{\cal P}_2}=\{A_{t,\tau,a,p_0}(\zeta); (t,\tau,a)\in {\cal P}_2,
p_0\in {\cal K} \cap N, \zeta \in \stackrel{\circ}{\Delta}_1\}$$ and
similarly for $N_{{\cal P}_2'}$.  Since the rank of the mapping
$${\cal P} \times ({\cal K}\backslash N)\times
\stackrel{\circ}{\Delta}_1 \ni (t, \tau, a, p_0, \zeta) \mapsto A_{t,
\tau, a, p_0}(\zeta) \in \C^n \backslash M$$ is equal to $(2n-2)$,
$N_{{\cal P}_2}$ and $N_{{\cal P}_2'}$ are closed two-codimensional
submanifolds of ${\cal W}_2, {\cal W}_2'$ respectively.  This implies
that
$$({\cal W}\backslash N_{{\cal P}_2}) \cap ({\cal W}\backslash
N_{{\cal P}_2'})= {\cal W}\backslash (N_{{\cal P}_2} \cup N_{{\cal
P}_2'})$$ is an open connected set. The restrictions of $F_2$ and
$F_2'$ to $({\cal W}\backslash N_{{\cal P}_2})$ and $({\cal
W}\backslash N_{{\cal P}_2'})$ therefore stick together in a function
$F$ that is holomorphic into
$${\cal W}\backslash (N_{{\cal P}_2} \cap N_{{\cal P}_2'}).$$ Indeed,
both are continuous up to $M\backslash N$, which is a uniqueness set,
and coincide there with $f|_{M\backslash N}$ by construction. Assuming
again that $\omega$ is thin near $0$, we insure that
$$\omega\cap({\cal W}\backslash (N_{{\cal P}_2} \cap N_{{\cal
P}_2'}))$$ is a connected open set. Therefore $f$ and $F$ stick
together in a well-defined holomorphic function into $\omega \cup
({\cal W}\backslash (N_{{\cal P}_2} \cap N_{{\cal P}_2'})).$

According to condition $(ii)$,
$$\overline{N}_{{\cal P}_2} \cap \overline{N}_{{\cal P}_2'} \cap M=
N.$$ Furthermore, let $\tilde{z}=A_{\tilde{t}, \tilde{\tau},
\tilde{a}, \tilde{p}_0}(\tilde{\zeta})$ be a point in $N_{{\cal
P}_2}\cap N_{{\cal P}_2'} \cap {\cal W}$ (if nonempty). Then,
according to conditions $(i)$ and $(ii)$, $T_{\tilde{z}}N_{{\cal
P}_2}$ and $T_{\tilde{z}}N_{{\cal P}_2'}$ intersect transversally in
$T_{\tilde{z}}\C^n$, since these are close to $T_0N + JT_0
\overline{\Gamma}_2$ and $T_0N+JT_0\overline{\Gamma}_2'$
respectively. Therefore, in a neighborhood $\tilde{\cal W}$ of
$\tilde{z}$ in ${\cal W}$, $L_{\tilde{z}}=N_{{\cal P}_2}\cap N_{{\cal
P}_2'}\cap \tilde{\cal W}$ consists of a four-codimensional manifold
in $\tilde{\cal W}$.  $L_{\tilde{z}}$ is removable for functions which
are holomorphic into $\tilde{\cal W} \backslash L_{\tilde{z}}$.  We
therefore showed that $N_{{\cal P}_2} \cap N_{{\cal P}_2'}$ is in fact
removable for functions holomorphic into ${\cal W}\backslash (N_{{\cal
P}_2}\cap N_{{\cal P}_2'})$. So $F$ extends holomorphically through
$N_{{\cal P}_2}\cap N_{{\cal P}_2'}$ as a function $F\in {\cal
O}({\cal W})$ with $F=f$ in the intersection of ${\cal W}$ with a
neighborhood of $M\backslash N$ in $\C^n$.

The proof of Proposition 4.1 is complete.

\smallskip
To achieve the proof of Theorem 1, we shall remark that in the proof
of Proposition 4.1, we only used the part of $\omega$ lying near
$A(-1)$ to perform deformations of $A$.  Therefore

\smallskip
\noindent
{\sc Proposition 4.2.}  {\it Let $M$ be generic,
$C^{2,\alpha}$-smooth, let $z_0\in M$, let $N\ni z_0$ be a $C^1$
submanifold with $codim_M N=3$ and $T_{z_0}N\not\supset T_{z_0}^cM$.
Assume there exists a sufficiently small embedded analytic disc $A\in
C^{2,\beta}(\overline{\Delta})$ attached to $M$, $A(1)=z_0$, with $A(b
\Delta \backslash \{1\}) \subset M \backslash N$,
$\frac{d}{d\theta}|_{\theta=0} A(e^{i\theta})=v_0\not\in T_{z_0}^cM$,
$v_0 \not\in T_{z_0}N$ and all discs in ${\cal V}(A,\delta)$ attached
to $M\backslash N$ are analytically isotopic to a point in
$M\backslash N$, for some $\delta>0$ and let ${\cal W}_0$ be an open
wedge attached to $M\backslash N$ containing a neighborhood of $A(-1)$
in $\C^n$.  Then there exists a wedge ${\cal W}$ of edge $M$ at $(z_0,
Jv_0)$ such that for every continuous function $f$ in ${\cal W}_0
\cup(M\backslash N)$ which is $CR$ on $M\backslash N$ and holomorphic
into ${\cal W}_0$, there exists a holomorphic function $F$ in ${\cal
O}({\cal W})$ continuous up to $M\backslash N$ with $F=f$ in a
neighborhood of $z_0$ in $M\backslash N$.}

\smallskip
{\it Proof.} We only have to check that the extension $F(z)$ defined
unambiguously in terms of a Cauchy integral on boundaries of discs
becomes holomorphic into the wedge ${\cal W}$ foliated by images of
the interior of $\Delta_1$.  One can introduce very small deformations
$M^d$ of $M$ into ${\cal W}_0$ depending on a real parameter $d\geq 0$
with $M^0=M$, apply the continuity principle argument to conclude
$F^d$ is holomorphic and let $d$ tend to $0$.

\smallskip
The proof of Theorem 1 can be completed as follows.  Proposition 2.1
and Proposition 3.7 insure the existence of an embedded disc $A$
meeting $N$ only at one point $A(1)$ of its boundary with nearby discs
being analytically isotopic to a point in $M\backslash N$.  Slightly
deform $M$ near $A(-1)$ in a manifold $M^d$ by pushing it into the
open wedge ${\cal W}_0$ of automatic extension, the deformation $A^d$
of $A$ remaining a disc with all similar properties.  Then Proposition
4.2 yields a holomorphic extension in a wedge of edge $M^d$ at $z_0$.
Since $M^d\equiv M$ near $z_0$, this gives a wedge of edge $M$ at
$z_0$.

The proof of Theorem 1 is complete.

\medskip
\noindent
{\sc Corollary 4.3.} {\it Let $M$ be a real analytic generic manifold
in $\C^n$ $(n\geq 2)$ of finite type at every point with $\hbox{CRdim}
\ M = p \geq 1$. Then every real analytic subset $A\subset M$ with
$\hbox{codim}_M A \geq 3$ is removable.}

\smallskip
{\it Proof.}  Let $A\subset M$ be a real analytic subset with
$\hbox{dim} \ A <2p$ and let $M^d\subset \C^n$ be a $C^{\infty}$
generic manifold of finite type at every point with $M^d\supset A$ and
$TM^d|_A=TM|_A$.  Fix an open neighborhood $V$ of $A$ in $M^d$ and
$\varepsilon >0$.  We shall show that there exists a $C^{\infty}$
generic manifold $M^{d_m}$ with $M^{d_m}\cap (M^d\backslash
V)=M^d\backslash V$ and $||M^{d_m}-M^d||_{C^{\infty}}<\varepsilon$
such that for each function $f\in \hbox{CR}(M^d\backslash N)$, there
exists a function $f^{d_m}\in \hbox{CR}(M^{d_m})$ with $f^{d_m}=f$ on
$M^d\backslash V$.  According to the Deformation Lemma below, this
will give the desired result on removability of $A$.

$A$ admits a stratification with the property that the closure of each
stratum only intersects strata of smaller dimension. Take $N'$ a
connected stratum of maximal possible dimension. Then $N'$ is an
embedded submanifold of $M^d\backslash (A\backslash N')$ to which
Theorem 1 applies, since $\hbox{dim} \ N' < 2p$ implies $T_zN'
\not\supset T_z^cM^d$ for every $z\in N'$. Since $N'$ is removable,
there exists a $C^{\infty}$ deformation $M^{d_1}$ over $N'$ of $M^d$
with support in $V$ such that
$||M^{d_1}-M^d||_{C^{\infty}}<\varepsilon_1<\varepsilon/2$ and
$M^{d_1}\cap (M^d\backslash V)=M^d \backslash V$ and $f^{d_1}\in
\hbox{CR}(M^{d_1} \backslash (A\backslash N'))$ with $f^{d_1}=f$ in
$M^d\backslash V$.  For $\varepsilon_1$ small enough, $M^{d_1}$ is of
finite type at every point.  Since the stratum of maximal possible
dimension in the set of remaining strata always looks like a locally
embedded submanifold, all strata of $A$ can be successively removed on
the successive deformations of $M^d$. Hence $A$ is removable.

Assume by induction that given a real analytic subset $A' \subset M$
with $\hbox{dim} \ A' \leq k \leq \hbox{dim} \ M-4$ and a $C^{\infty}$
everywhere of finite type generic manifold in $\C^n$, $M^d \supset
A'$, with $TM^d|_{A'}=TM|_{A'}$, $A'$ is removable for CR functions on
$M^d\backslash A'$.  Let $A\subset M$ be a stratified real analytic
subset, $\hbox{dim} \ A=k+1$ and let $M^d$ be a $C^{\infty}$
everywhere of finite type generic manifold in $\C^n$ with $M^d\supset
A$ and $TM^d|_A=TM|_A$. Let $V$ be a neighborhood of $A$ in $M^d$ and
let $\varepsilon>0$ be arbitrary. Choose $N$ a connected stratum of
$A$ of maximal dimension.  Then the set
$$N^c=\{z\in N; \ T_zN \supset T_z^cM^d=T^c_zM\}$$ is a proper
subanalytic set of $N$. Indeed, since $M^d$ is minimal at every point,
$M^d$ does not contain germs of CR manifolds $\Sigma$ with
$\hbox{CRdim} \ \Sigma =\hbox{CRdim} \ M^d$.  A relatively open set in
$N$ contained in $N^c$ would be such a CR manifold by definition of
$N^c$.  Therefore, the dense open subset of $N$
$$N\backslash N^c=\{z\in N; \ T_zN \not\supset T_z^cM \}$$ is
removable, by virtue of Theorem 1.

Fix a function $f\in \hbox{CR}(M^d\backslash N)$.  Smoothly deform
$M^d$ in a manifold $M^{d_1}$ of class $C^{\infty}$ and of finite type
at every point over $N\backslash N^c$ by pushing it into the wedge
where removability of $N\backslash N^c$ holds and keeping both
$TM^d|_{N^c}=TM^{d_1}|_{N^c}$ and $M^{d_1} \cap (M^d\backslash V) =
(M^d\backslash V)$, with $||M^{d_1} - M^d||_{C^{\infty}}
<\varepsilon_1<\varepsilon/2$.  Take the restriction $f^{d_1}$ of $f$
to $M^{d_1}\backslash N^c$.  According to the induction hypothesis,
$N^c$ is removable on $M^{d_1}$, hence there exists a smooth
deformation $M^{d_2}$ of $M^{d_1}$ with
$||M^{d_2}-M^{d_1}||_{C^{\infty}}<\varepsilon_2<\varepsilon/4$,
$TM^{d_2}|_{A\backslash N}=TM^d|_{A\backslash N}$, $M^{d_2}\cap
(M^d\backslash V)=M^d\backslash V$ such that for each function $f\in
\hbox{CR}(M^d\backslash A)$ there exists $f^{d_2}\in
\hbox{CR}(M^{d_2}\backslash (A\backslash N))$ with $f^{d_2}=f$ on
$M^d\backslash V$.

All strata $N$ can be successively removed by applying the induction
hypothesis to $N^c$.

The proof finishes out with the following.

\smallskip
\noindent
{\sc Deformation Lemma.} {\it Let $M$ be minimal at $z_0\in M$, $\Phi$
a proper closed subset with $b\Phi\ni z_0$. Assume that for each
neighborhood $V$ of $\Phi$ in $M$, each $\varepsilon>0$, there exists
a manifold $M^d$ with $M^d\cap (M\backslash V)=M\backslash V$,
$||M^d-M||_{C^{2,\alpha}}<\varepsilon$ such that for each function
$f\in\hbox{CR}(M\backslash \Phi)$ there exists a function
$f^d\in\hbox{CR}(M^d)$ with $f^d\equiv f$ on $M\backslash V$. Then
there exists a wedge ${\cal W}$ of edge $M$ at $z_0$ with
$\hbox{CR}(M\backslash \Phi)$ extending holomorphically into ${\cal
W}$ and continuously in ${\cal W} \cup (M\backslash \Phi)$.}

\smallskip
{\it Proof.} The result follows from the existence of a disc with
defect $0$ attached to $M^d$ when $\varepsilon$ is small enough,
obtained as a perturbation of a disc with defect $0$ attached to the
manifold $M$, minimal at $z_0$.  Each $f^d$ extends holomorphically
into a wedge ${\cal W}^d$ of edge $M^d$ which stabilizes as
$\varepsilon$ tends to $0$ and gives a wedge ${\cal W}$ attached to
$M$ as $V$ shrinks to $\Phi$.

\smallskip
The proof of Corollary 4.3 is complete.

\smallskip
\noindent
{\bf 5. Proof of Theorem 2.} The meaning of Theorem 2 is that, when
the Cauchy-Riemann dimension of $M$ is greater or equal to $2$, some
two-codimensional singularities are locally removable.  However, this
is not generically true in CR dimension equal to $1$. Take for example
$M$, a hypersurface in $\C^2$ and $g$ a holomorphic function in a
neighborhood of some $z_0\in M$ in $\C^2$ such that $g(z_0)=0$ and
$dg(z_0)\neq 0$. Whenever $\Sigma_g=\{g=0\}$ is not tangent to $M$ at
$z_0$, $\frac{1}{g}|_{M\backslash N}$ defines a nonextendible CR
function in any side, and $N=\Sigma_g\cap M$ is locally a
two-codimensional submanifold in $M$ with $T_{z_0}N\cap
T_{z_0}^cM=\{0\}$. The difference between CR dimension one and CR
dimension greater than two is explained in the two proofs of the
following theorem.

\smallskip
\noindent
{\sc Theorem 5.A.} {\it Let $M$ be a $C^{2,\alpha}$-smooth
$(0<\alpha<1)$ generic manifold in $\C^n (n\geq 3)$ with $\hbox{CRdim}
\ M=p \geq 2$ and $N \subset M$ a $C^{2,\alpha}$ submanifold of $M$
with $codim_M N =2$ that is generic in $\C^n$. Assume that there
exists an open neighborhood $V$ of $N$ in $M$ such that $(M\backslash
N)\cap V$ has the wedge extension property.  Then there exists an open
wedge ${\cal W}$ attached to $M \cap V$ such that every continuous CR
function on $(M\backslash N)\cap V$ extends holomorphically into
${\cal W}$ and continuously in $(M\backslash N) \cup {\cal W}$.}

\smallskip
{\it Proof.}  The condition that $N$ is a generic manifold in $\C^n$
means that $\hbox{CRdim} \ (T_n N \cap T_n^c M)=p -2, \forall n \in
N$.  Therefore, $N$ cannot contain germs of CR manifolds with CR
dimension equal to $p -1$. In two senses of speaking, Theorem 5.A
treats the generic case of Theorem 2.
 
To obtain Theorem 5.A, we shall use a natural and beautiful
deformation result due to J\"oricke, obtained in \cite{JO3} as a tool
in her deriving a proof of a conjecture of Tr\'epreau.

\smallskip
\noindent
{\sc Theorem.} ({\sc J\"oricke} \cite{JO3}).  {\it Let $z_0 \in M$ be
generic, $C^{2, \alpha}$-smooth $(0<\alpha <1)$ in $\C^n$ and let $C$
be a truncated open convex cone in $\C^n$ with vertex $z_0$ and some
$v \in T_{z_0}^cM \backslash \{0\}$ in $C$.  Then for every
neighborhood $\omega \subset C$ of $C \cap M$ in $\C^n$, every $\beta,
0<\beta <\alpha,$ there exists a $C^{2, \beta}$-smooth generic
manifold $M^d \subset M \cup \omega$ with $M^d \backslash
\omega=M\backslash \omega$ such that $M^d$ is minimal at $z_0$.}

\smallskip
Let $z_0\in N$ and let $U\subset V$ be a small neighborhood of $z_0$
in $M$.  Let $M_1$ be a closed $C^{2,\alpha}$ one codimensional
connected generic submanifold of $M\cap U$ containing $N\cap U$. Since
$N$ itself is generic, there exists $C$ a truncated open convex cone
in $\C^n$ with vertex $z_0$ and some $v\in T_{z_0}^cM_1\backslash
\{0\}$ in $C$ such that $C\cap N=\emptyset$.

Since $(M\cap U)\backslash N$ has the wedge extension property, there
exists a wedge ${\cal W}_0$ attached to $(M\cap U)\backslash N$ to
which CR functions on $M\backslash N$ holomorphically extend. We can
first perform a small $C^{2,\alpha}$ deformation $M^d$ of $M$ into
${\cal W}_0$ leaving $M \backslash C$ fixed in order that ${\cal W}_0$
becomes a neighborhood $\omega$ of $(M_1)^d\cap C$ in
$\C^n$. Moreover, we can assume that the tangent spaces and the
complex tangent spaces to $M$ and $M_1$ at $z_0$ are fixed under
$d$. Thus, $v\in T_{z_0}^c (M_1)^d$ too.

Secondly, minimalize $(M_1)^d$ at $z_0$ by applying the deformation
theorem of J\"oricke: we get a manifold $(M_1)^{d_2}$ contained in
$\omega\cup (M_1\backslash C)$ such that $(M_1)^{d_2}\cap
(M_1\backslash C)=M_1\backslash C$ and $(M_1)^{d_2}$ is minimal at
$z_0$. Furthermore, this deformation can be extended in a smooth
deformation $M^{d_2}$ of $M$ with support in $\overline{C\cap M}$,
{\em i.e.} $M^{d_2}\cap (M\backslash C)=M\backslash C$.

We shall now show that there exists a neighborhood $W$ of $z_0$ in $M$
such that every point in $(W\cap N)\backslash \{z_0\}$ is removable,
that is, for each point $z\in (W\cap N)\backslash \{z_0\}$ there
exists a wedge ${\cal W}_z$ of edge $M$ at $z$ with $CR((M\backslash
N)\cap V)$ extending holomorphically into ${\cal W}_z$.

Indeed, let $f$ be a function in CR$((M\backslash N)\cap V)$.  Since
$f$ extends holomorphically into ${\cal W}_0$, we get a function,
still denoted by $f$, which is CR on $M^{d_2}\backslash N$.  We shall
show in Lemma 5 below that $f$ extends holomorphically into a wedge of
edge $M^{d_2}$ at $z_0$. Since $M\equiv M^{d_2}$ in a neighborhood of
each point $z\in N$ with $z\neq z_0$, this gives the desired result.

\smallskip
\noindent
{\sc Lemma.} ({\sc Chirka} and {\sc Stout} \cite{CS}).  {\it Let $M$
be generic, $C^{2,\alpha}$ and let $M_1\subset M$ be a
$C^{2,\beta}$-smooth one codimensional submanifold, generic and
minimal at $z_0\in M_1$ in $\C^n$. Let $\Phi \subset M_1$ be a proper
closed subset of $M_1$, $b\Phi \ni z_0$.  Then $\Phi$ is removable at
$z_0$.}

\smallskip
{\it Proof.}  Fix a function $f\in \hbox{CR}(M\backslash \Phi)$.
Include $M_1$ in a regular one parameter family $M_{1,s}, |s|<\delta,
\delta >0$ of $C^{2, \beta}$ manifolds contained in $M$, such that
$M_{1,0} = M_1,$ $M_{1,s} \cap \Phi = \emptyset, s \neq 0$ and $\cup_s
M_{1,s}= M$ near $z_0$. Using cartesian equations for both $M_1$ and
$M$, these manifolds $M_{1,s}$ can be defined as translation-like
modifications of $M_1$.  Furthermore, we can impose that all manifolds
$M_{1,s}$ contain some fixed point $z_1\in M_1\backslash \Phi$.  The
method of sweeping out by wedges can be applied (\cite{CS}).  Since
existence of a disc with minimal defect is a stable property under
$C^{2, \beta}$ deformations, Tumanov's theorem gives: every $f_s =
f|_{M_{1,s}}$ is wedge extendible into a wedge ${\cal W}_s$ of edge
$M_{1, s}$ for $s \in (-\delta ', \delta') \backslash \{ 0 \}$, $0 <
\delta' \leq \delta$, and the wedges depend smoothly on $s$.  By the
uniqueness theorem and the fact that $M\backslash \Phi$ is connected
near $z_0$, all $M_{1,s}$ containing $z_1\in M_1\backslash \Phi$,
holomorphic extensions obtained constitute a single function
$\tilde{f}$ holomorphic into the union ${\cal W} = \cup_{s \neq 0}
{\cal W}_t$.  By construction, the set ${\cal W}$ contains a
nontrivial wedge of edge $M$ at $z_0$.  Moreover, $\tilde{f}$ is
continuous up to $M\backslash \Phi$ and equals $f$ there.

The proof of Theorem 5.A is complete.

\smallskip
{\it Remark.} The method of sweeping out by wedges can only be applied
when $\hbox{CRdim} \ M \geq 2$, since otherwise any generic manifold
$M_1$ as above is totally real in $\C^n$ hence has no CR structure.

\smallskip
{\it Remark.} The method of proof of Theorem 5.A cannot be used to
derive Theorem 1 or Theorem 2 in full generality for the following
reason.  If $N$ contains a CR manifold $\Sigma$ through $z_0$ with
$\hbox{CRdim} \ \Sigma=p-1$, no one codimensional submanifold $M_1$ of
$M$ containing $N$ that is generic in $\C^n$ can be minimal, since
then $\Sigma$ contains the local CR orbit ${\cal
O}^{loc}_{CR}(z_0,M_1)$.  Notice that $\Sigma$ can be a proper
submanifold of $N$, $2p-2\leq dim \ \Sigma \leq dim \ N=2p+q-2$ or
$2p+q-3$.

\smallskip
A second proof of Theorem 5.A can be obtained as follows and allows
$N$ to be of class $C^1$.

\smallskip
\noindent
{\sc Proposition 5.B.} {\it Let $M$ be a $C^2$-smooth generic manifold
in $\C^n$, $\hbox{CRdim} \ M=p\geq 2$, let $N \subset M$ be a
$C^1$-smooth {\rm generic} submanifold with $\hbox{codim}_M N = 2$ and
let $z_0\in N$.  Then there exist two neighborhoods $U,V \subset
\subset U$ of $z_0$ in $M$ such that each function $f \in
\hbox{CR}(U\backslash N)$ can be uniformly approximated on compact
subsets of $V\backslash N$ by holomorphic polynomials.}

\smallskip
{\it Proof.} This is an adaptation of the approximation theorem of
Baouendi and Treves. Let $L_0$ be a maximally real submanifold through
$z_0$ with $T_{z_0}L_0 = \R^n$ contained in $N$, let $H$ be a $C^2$
manifold through $z_0$ with $\hbox{dim} \ H + \hbox{dim} \ L_0 =
\hbox{dim} \ M$ and $T_{z_0} L_0 + T_{z_0} H = T_{z_0} M$. There are
maximally real manifolds $L_h$ through $h \in H$, closed in a fixed
neighborhood $U$ of $z_0$ in $M$, with the properties that $M\cap V$
is the disjoint union of the $L_h \cap V$, for some neighborhood $V$
of $z_0$ in $M$, the $L_h$ are uniformly close in $C^1$ norm to $L_0$,
$L_h \subset N$ if $h\in N$ and $L_h \cap N = \emptyset$ if $h \not
\in N$.  Fix a manifold $L_1$ contained in $M\backslash N$ and set, if
$f\in \hbox{CR}(M\backslash N)$ and $\widehat{z}\in V\backslash N$
$$G_{\tau} f( \widehat{z}) = \left(\frac{\tau}{\pi}\right)^{n/2}
\int_{L_1} e^{-\tau (z-\widehat{z})^2} f(z)dz,$$ where
$(z-\widehat{z})^2=(z_1-\widehat{z}_1)^2+\cdots+
(z_n-\widehat{z}_n)^2$, $\tau>0$ and $dz= dz_1 \wedge \cdots \wedge
dz_n$.

If $\widehat{z} \in V\backslash N$, $\widehat{z}$ belongs to some
manifold $L_{\hat{h}}$ and we have, if $f$ is of class $C^1$
$$G_{\tau} f( \widehat{z}) = \left(\frac{\tau}{\pi}\right)^{n/2}
\int_{L_{\hat{h}}} e^{-\tau (z-\widehat{z})^2}
f(z)dz+\left(\frac{\tau}{\pi}\right)^{n/2} \int_{\Sigma}
d(e^{-\tau(z-\widehat{z})^2} f(z)dz) =
\left(\frac{\tau}{\pi}\right)^{n/2} \int_{L_{\hat{h}}} e^{-\tau
(z-\widehat{z})^2} f(z)dz,$$ by Stokes' theorem. The last equality
holds since $f$ is CR on $M\backslash N$, hence $d(f(z)dz)=0$, and
since $\Sigma$ with $b \Sigma=L_1-L_{\hat{h}}$ can be chosen to be
contained in $M\backslash N$. When $f$ is only continous, the middle
equality has to be interpreted in the distribution sense.

Analysing the real and imaginary parts of the phase function
$-\tau(z-\widehat{z})^2$ on $L_{\hat{h}}$, one can show that the last
integral tends to $f(\widehat{z})$ as $\tau$ tends to $\infty$ if the
$L_{\hat{h}}$ are sufficiently close to $L_0$ in $C^1$ norm. Then the
expression defining $G_{\tau} f( \widehat{z})$ by integration on $L_1$
gives converging polynomial sequences using a truncated developement
of the exponential in power series.

The proof of Proposition 5.B is complete.

\smallskip
{\it Remark.} The approximation property is still valid for CR
functions extending holomorphically into some wedge ${\cal W}_1$ of
edge $V$ at $z_1\in V$, on compact subsets of $(V\backslash N)\cup
{\cal W}_2, {\cal W}_2\subset {\cal W}_1$. Indeed, some uniformity is
allowed in choosing the $L_h$ filling $(V\backslash N)\cup {\cal
W}_2$.

\smallskip
Here is a second version of Theorem 5.A, whose proof uses deformations
of discs instead of the minimalization theorem and the sweeping out by
wedges lemma.

\smallskip
\noindent
{\sc Theorem 5.A.1} {\it Let $M$ be a $C^{2,\alpha}$-smooth
$(0<\alpha<1)$ generic manifold in $\C^n (n\geq 3)$ with $\hbox{CRdim}
\ M=p \geq 2$ and $N \subset M$ a $C^1$ submanifold of $M$ with
$codim_M N =2$ that is generic in $\C^n$. Assume there exists an open
neighborhood $V$ of $N$ in $M$ such that $(M\backslash N)\cap V$ has
the wedge extension property.  Then there exists an open wedge ${\cal
W}$ attached to $M \cap V$ such that every continuous CR function on
$(M\backslash N)\cap V$ extends holomorphically into ${\cal W}$ and
continuously in $(M\backslash N) \cup {\cal W}$.}

\smallskip
{\it Proof.} According to Proposition 5.B, CR functions on
$M\backslash N$ are locally uniformly approximable on compact subsets
of $M\backslash N$ by holomorphic polynomials.  However, embedded
discs attached to $M\backslash N$ are not in general analytically
isotopic to a point in $M\backslash N$, since the first homology
groups $H_1(V\backslash N)\equiv \Z$, for $V$ a small open ball of
center $z_0 \in N$.

Theorem 5.A.1 will be a consequence of Proposition below and existence
of a good disc.  Existence is checked as follows.  There are
holomorphic coordinates $(w,z)$ on $\C^n$ as in $(10)$ such that $T_0
N =\{v_1=v_2=0\}$. Then for small $c>0$, the disc $A_c$ with
holomorphic $w$-component $w_c(\zeta)=(c(1-\zeta),
ic(1-\zeta),0,...,0)$ and $z$-component satisfying the Bishop equation
$x_c=-T_1h(w_c, x_c)$ is attached to $M$ and satisfies $A_c(1)=0\in N,
A_c(b\Delta \backslash \{1\}) \subset M\backslash N, v_0=
\frac{d}{d\theta}|_{\theta=0} A_c(e^{i\theta}) \not\in T_0N$.

\smallskip
\noindent
{\sc Proposition 5.1} {\it Let $M$ be generic, $C^{2, \alpha}$-smooth,
$p\geq 1$, let $z_0 \in M$, let $N\ni z_0$ be a $C^1$ submanifold with
$\hbox{codim}_M N =2$ and $T_{z_0} N \not\supset T_{z_0}^c M$ and let
$V$ be a neighborhood of $z_0$ in $M.$ Assume that there exists a
sufficiently small analytic disc $A\in C^{2, \beta}
(\overline{\Delta})$ attached to $M$ with $A(1)=z_0, A(b
\Delta\backslash \{1\}) \subset V\backslash N,
\frac{d}{d\theta}|_{\theta=0} A(e^{i\theta})=v_0 \not\in T_{z_0}^cM,
v_0 \not\in T_{z_0}N$.  Let ${\cal W}_1$ be a wedge of edge $M$ at
$A(-1)$ and assume that CR functions on $M\backslash N$ extending
holomorphically into ${\cal W}_1$ are uniformly approximable on
compact subsets of $(V\backslash N)\cup {\cal W}_2, {\cal W}_2\subset
{\cal W}_1$, by holomorphic polynomials.  Then for each $\varepsilon
>0$, there exists $v_{00}\in T_{z_0}M$ with
$|v_{00}-v_0|<\varepsilon$, $v_{00}\not\in T_{z_0}N$, $v_{00}\not\in
T_{z_0}^c M$ and a wedge ${\cal W}$ of edge $M$ at $(z_0, Jv_{00})$
such that, if a function $f\in \hbox{CR}(M\backslash N)$ extends
holomorphically into ${\cal W}_1$, it extends to be holomorphic into
${\cal W}$.}

\smallskip
{\it Remark.} The proposition holds even if $N$ is not generic. Thus,
in a generic manifold with $M\backslash N$ having the wedge extension
property, a $C^1$ two-codimensional singularity $N$ with $T_{z_0} N
\not\supset T_{z_0}^cM$ is removable at $z_0$ {\it if and only if} CR
functions on $(M\backslash N)\cup {\cal W}_1$ are locally uniformly
approximable by holomorphic polynomials on compact subsets of
$(M\backslash N)\cup {\cal W}_2$ near $z_0$.

\smallskip
{\it Proof.} The proof uses the same deformations of discs as in
Section 4, but holomorphic extendibility into the sets ${\cal
W}\backslash N_{\cal P}$ now is a direct consequence of the uniform
approximability of $\hbox{CR}(M\backslash N)$ by holomorphic
polynomials, and the maximum principle, as usual in the field
\cite{TR1} \cite{TU1}.
 
If $v_0\in T_{z_0}^cM$, choose a disc $A_{00}=A_{t_1,\tau_1,a_1}$ in
the family of deformed discs constructed in $(22)$ with
$v_{00}=\frac{d}{d\theta}|_{\theta=0}
A_{00}(e^{i\theta})=v_{00}\not\in T_{z_0}^cM$, $v_{00}\not\in
T_{z_0}N$ and $|v_{00}-v_0|<\varepsilon.$
 
The argument of overlapping wedges gives the following.  Condition
$(i)$ insures that ${\cal W}_2$ is a wedge of edge $M$ at $(0,
Jv_{00})$ and
$$N_{{\cal P}_2}=\{A_{t,\tau,a, p_0}(\stackrel{\circ}{\Delta}_1); \
(t,\tau,a,p_0)\in {\cal P}_2, p_0 \in {\cal K} \cap N \}$$ is a closed
one-codimensional conic submanifold in ${\cal W}_2$ and similar
properties hold for ${\cal W}_2'$, $N_{{\cal P}_2}'$.  Connectedness
arguments run without modifications, since $Jv_{00} \not\in T_0M$.

Given a function $f\in \hbox{CR}(M\backslash N)$, one obtains a
function $F$ that is holomorphic into
$$({\cal W}_2 \cup {\cal W}_2')\backslash (N_{{\cal P}_2} \cap
N_{{\cal P}_2}')$$ and continuous up to $M\backslash N$, with $F=f$
near $0$ on $M\backslash N$.

Condition $(i)$ insures that the two wedges ${\cal W}_2$ and ${\cal
W}_2'$ contain a wedge ${\cal W}$ of edge $M$ at
$(0,Jv_{00})$. Condition $(ii)$ together with the facts that $T_0 N
\not\supset T_0^cM$ and $v_{00} \not\in T_0 N$ imply that
$M_1=N_{{\cal P}_2} \cap N_{{\cal P}_2}'$ is locally ({\em i.e.} near
a point in ${\cal W}$, see the proof of Proposition 4.1) a $C^1$
two-codimensional submanifold which is generic in $\C^n$.

A generic locally closed two -codimensional manifold $M_1$ in an open
set ${\cal W}$ is removable for functions which are holomorphic into
${\cal W}\backslash M_1$.

Indeed, according to Proposition 5.B, CR ({\em i.e. holomorphic})
functions on ${\cal W}\backslash M_1$ are locally uniformly
approximable by holomorphic polynomials.  Therefore, given a function
$f\in {\cal O}({\cal W}\backslash M_1)$ and a small disc $A$ attached
to ${\cal W}\backslash M_1$, $f\circ A|_{b\Delta}$ extends
holomorphically into $\Delta$. Let $z_0\in M_1\cap {\cal W}$ and
choose a holomorphic coordinate system at $z_0$ such that $T_{z_0}
M_1=\{y_1=y_2=0\}$. Then the discs $A_{c, p_0}=(c\zeta, ic\zeta,
0,...,0) +p_0, c>0, p_0 \in M_1, |p_0| <<c$ satisfy $A_{c,
p_0}(1)=p_0\in M_1, A_{c,p_0}(b\Delta) \subset {\cal W}\backslash
M_1$.  As a consequence of the existence of such a family, each
function $f\in {\cal O}({\cal W}\backslash M_1)$ extends continuously,
hence holomorphically, through $M_1$, as a $F\in{\cal O}({\cal W})$
with $F=f$ on ${\cal W}\backslash M_1$, and with
$$F(p_0)=\frac{1}{2i\pi} \int_{b\Delta} \frac{f\circ A_{c,
p_0}(\eta)}{\eta-\zeta} d\eta$$ for $p_0 \in M_1\cap {\cal W}$.

 The proof of Proposition 5.1 is complete.

\smallskip
Return to the proof of Theorem 2. Since $N$ does not consist of a CR
manifold with $\hbox{CRdim} \ N=p-1$, all points where $N$ is generic
in $\C^n$ are removable.  Therefore, the set
$$N^{CR}=\{n\in N; \ \hbox{CRdim} \ T_n N=p-1 \}$$ is a proper closed
subset of $N$.  The following Theorem finishes out the proof of
Theorem 2.

\smallskip
\noindent
{\sc Theorem 5.2.}  {\it Let $M$ be a $C^{2,\alpha}$-smooth
$(0<\alpha<1)$ generic manifold in $\C^n (n\geq 3)$ with
$p=\hbox{CRdim} \ M \geq 2$, $N \subset M$ a connected $C^1$
submanifold of $M$ with $codim_M N =2$ such that $T_z N \not\supset
T_z^cM$ for each $z\in N$ and let $\Phi$ be a proper closed subset of
$N$. Assume that there exists a neighborhood $V$ of $N$ in $M$ such
that $M$ is minimal at every point of $V$.  Then there exists an open
wedge ${\cal W}$ attached to $M \cap V$ such that every continuous CR
function on $(M\backslash \Phi)\cap V$ extends holomorphically into
${\cal W}$ and continuously in $(M\backslash \Phi) \cup {\cal W}$.}

\smallskip
{\it Proof.}  Fix a function $f\in CR(M\backslash \Phi)$.
$(M\backslash \Phi)\cap V$ has the wedge extension property, since $M$
is minimal at every point of $V$.  Furthermore, $f$ extends
holomorphically into some open wedge ${\cal W}_0$ attached to
$(M\backslash \Phi_{nr})\cap V$, where $\Phi_{r}$ denotes the set of
removable points of $\Phi$ and $\Phi_{nr}=\Phi\backslash \Phi_r$
denotes the set of nonremovable points.

According to Theorem 5.A.1, $\Phi_r$ contains all points of $\Phi$
where $N$ is generic in $\C^n$.  By definition, $\Phi_r$ is a
relatively open subset of $\Phi$ and $\Phi_{nr}$ is a proper closed
subset of $N$. Assume that $\Phi_{nr}$ is nonempty. We shall reach a
contradiction.

Since $N$ is connected, there exists a point $z_0$ in the relative
boundary of $\Phi_{nr}$ with respect to $N$.

Smoothly deform $M$ in a $C^{2,\alpha}$-smooth manifold $M^d$ by
pushing it slightly into ${\cal W}_0$, the deformation being
sufficiently small in order that there still exists a good disc
attached to $M^d$ at $z_0$, as constructed in Section 2 on everywhere
minimal CR manifolds.

According to the Deformation Lemma, Proposition 2.1 and Proposition
3.8, removability of $z_0$ will be a consequence of the following.

\smallskip
\noindent
{\sc Proposition 5.3.} {\it Let $M$ be generic, $C^{2,\alpha}$-smooth,
let $z_0\in M$, let $N\ni z_0$ be a $C^1$ submanifold with $codim_M
N=2$ and $T_{z_0}N\not\supset T_{z_0}^cM$ and let $\omega$ be a
neighborhood of $M \backslash \Phi$ in $\C^n$ for some proper closed
subset $\Phi$ of $N$ with $z_0\in b\Phi$.  Assume there exists a
sufficiently small embedded analytic disc $A\in
C^{2,\beta}(\overline{\Delta})$ attached to $M$, $A(1)=z_0$, with $A(b
\Delta \backslash \{1\}) \subset M \backslash N$,
$\frac{d}{d\theta}|_{\theta=0} A(e^{i\theta})=v_0\not\in T_{z_0}^cM$,
$v_0 \not\in T_{z_0}N$ and all discs in ${\cal V}(A,\delta)$ attached
to $M\backslash \Phi$ are analytically isotopic to a point in
$M\backslash \Phi$, for some $\delta>0$.  Then there exists a wedge
${\cal W}$ of edge $M$ at $(z_0, Jv_0)$ such that for every
holomorphic function $f \in {\cal O}(\omega)$ there exists a function
$F \in {\cal O}({\cal W})$ with $F=f$ in the intersection of ${\cal
W}$ with a neighborhood of $M\backslash N$ in $\C^n$.}

\smallskip
{\it Proof.} The proof uses same deformations of discs as in Section 4
until we reach the argument of overlapping wedges. We take the
notations of Section 4.

Notice that condition $(i)$ insures that ${\cal W}_2$ is a wedge of
edge $M$ at $(0, \eta_0), \eta_0=Jv_0 \ \hbox{mod} \ T_{0}M$ and
$$N_{{\cal P}_2}=\{A_{t,\tau,a,p_0}(\stackrel{\circ}{\Delta}_1); \
(t,\tau,a) \in {\cal P}_2, p_0 \in {\cal K} \cap N \}$$ is a closed
one-codimensional conic submanifold in ${\cal W}_2$ (in other words, a
CR-wedge over $N$ or a {\it manifold with edge} $N$), $\Phi_{{\cal
P}_2}= \{A_{t,\tau,a,p_0}(\Delta_1); \ (t,\tau,a) \in {\cal P}_2, p_0
\in {\cal K} \cap \Phi \}$ being a proper closed subset of $N_{{\cal
P}_2}$. Since ${\cal W}_2 \backslash \Phi_{{\cal P}_2}$ is therefore
connected, the continuity principle argument and Cauchy's integral
yield a function $F_2$ that is holomorphic into ${\cal W}_2 \backslash
\Phi_{{\cal P}_2}$, and similarly also $F_2' \in {\cal O}({\cal
W}_2'\backslash \Phi_{{\cal P}_2}')$.  By shrinking $\omega$ near $0$
and ${\cal P}_2, {\cal P}_2'$ if necessary, we can insure that all the
open sets \begin{equation} \omega \cap ({\cal W}_2 \backslash
\Phi_{{\cal P}_2}) \ \ \ \ \ \omega \cap ({\cal W}_2'\backslash
\Phi_{{\cal P}_2'}) \ \ \ \ \ \omega \cap ( ({\cal W}_2 \cup {\cal
W}_2')\backslash (\Phi_{{\cal P}_2} \cap \Phi_{{\cal P}_2'}))
\end{equation} are connected.  Indeed, for the first two, this is true
if $\omega$ is sufficiently thin near $0$ and for the third, condition
$(i)$ forces the two connected wedges ${\cal W}_2$ and ${\cal W}_2'$
to be with nonempty intersection.  Since $F_2$ and $F_2'$ by
construction assume the values of $f$ on $M\backslash \Phi$, and since
$M\backslash \Phi$ is a uniqueness set, we have shown that there
exists a function $F$ that is holomorphic into \begin{equation} \omega
\bigcup (({\cal W}_2 \cup {\cal W}_2' ) \backslash (\Phi_{{\cal P}_2}
\cap \Phi_{{\cal P}_2'})) \end{equation} and $F|_{\omega}=f$.

Condition $(i)$ insures that the two wedges ${\cal W}_2$ and ${\cal
W}_2'$ contain a wedge ${\cal W}$ of edge $M$ at $(0, \eta_0)$. Take
for convenience the restriction of $F$ to $\omega \cup {\cal W}
\backslash (\Phi_{{\cal P}_2} \cap \Phi_{{\cal P}_2'})$ and still
denote it by $F$.  Let $\tilde{z}=A_{\tilde{t},\tilde{\tau},\tilde{a},
\tilde{p}_0}(\tilde{\zeta})$ be a point in $\Phi_{{\cal P}_2}\cap
\Phi_{{\cal P}_2'}\cap {\cal W}$ (if such exists).  Let $\tilde{\cal
W}$ be a neighborhood of $\tilde{z}$ in ${\cal W}$.  Then, according
to conditions $(i)$ and $(ii)$ together with the facts that
$T_0N\not\supset T_0^cM$ and the projection on the $v_1$-axis of
$v_0=\frac{d}{d\theta} A_{0,0,0,0}(1)$ is nonzero, $N_{{\cal P}_2}
\cap N_{{\cal P}_2'}\cap \tilde{\cal W}$ is contained in a
two-codimensonal generic manifold $L_{\tilde{z}}$ passing through
$\tilde{z}$.  $L_{\tilde{z}}$ being generic, we can remove it for CR
functions which are holomorphic into $\tilde{\cal W}\backslash
L_{\tilde{z}}$. Indeed, this is done as in the proof of Proposition
5.1.  Hence also $\Phi_{{\cal P}_2}\cap \Phi_{{\cal P}_2'}\cap {\cal
W}$ is removable.  We therefore showed that $F$ extends
holomorphically through $\Phi_{{\cal P}_2}\cap \Phi_{{\cal P}_2'}$ as
a function $F\in {\cal O}({\cal W})$ continuous up to $M\backslash N$
with $F|_{M\backslash N}=f|_{M\backslash N}$.

The proof of Theorem 5.2 is complete.

\medskip
{\it Proof of Theorem 3.}  Theorem 3 is a corollary of the following.

\smallskip
\noindent
{\sc Theorem 5.4.} {\it Let $M$ be a $C^{2, \alpha}$-smooth
$(0<\alpha<1)$ generic manifold in $\C^n$ $(n\geq 3)$, let $N$ be a
closed connected $C^2$ generic submanifold of $M$ with $\hbox{codim}_M
N =1$ and let $\Phi\subset N$ be a proper closed subset of $N$. Assume
that $\hbox{CR}(M\backslash \Phi)$ has the wedge extension property
and let $\Phi_r$ denote the set of removable points in $\Phi$. Then
$\Phi_{nr}=\Phi\backslash \Phi_r$ has the following structure: $b
\Phi_{nr}$ is a union of CR orbits of $N$.}

\smallskip
{\it Proof.} Assume on the contrary that there exists $z_0\in b
\Phi_{nr}$ and a point $z_1\in {\cal O}_{CR}(N,z_0)$ with $z_1 \not\in
b \Phi_{nr}$. Either there exists such a $z_1$ with $z_1\in
M\backslash N$ or ${\cal O}_{CR}(N,z_0)\subset \Phi_{nr}$ and $z_1 \in
\hbox{Int}\ \Phi_{nr}$.  In the latter case, moving backwards along
some perturbation of the piecewise smooth integral curve of $T^c M$
joining $n$ with $z_1$, one obtains that there exists a $C^2$ integral
curve $s \mapsto X_s(z_2), s\in [0, s_0], s_0 >0$ of a $T^cM$-tangent
vector field $X$ with $z_2\in M\backslash N, X_s(z_2)\in M\backslash
N, s\in [0, s_0)$ and $z=X_{s_0}(z_2)\in b\Phi_{nr}$.  This is also
true in the first case.

Though $b\Phi_{nr}$ can make a too thin \'etroiture at $z$, there
exist points $z_0\in b \Phi_{nr}$ close to $z$ such that there exists
a truncated convex open cone $\Gamma$ at $(z_0, -X(z_0))$ contained in
$N$ with $\Gamma\cap \Phi_{nr}=\emptyset$.  Indeed, choose a real
euclidean coordinate system $(n_1,...,n_k), k=\hbox{dim }N$ on $N$
near $z$ such that integral curves of $X$ correspond to lines
$n_2=ct,...,n_k=ct.$ Let $s_1 < s_0$ be close to $s_0$, let $r>0$ be
so small that the closed ball $B_1=\{z\in N; \
(n_1-n_1^1)^2+\cdots+(n_k-n_k^1)^2 \leq r^2\}$ is contained in
$N\backslash \Phi$ where $z_1=X_{s_1}(z_2)$ has coordinates
$(n_1^1,...,n_k^1)$.  When $\delta \geq 1$ increases, the domains
$B_{\delta}=\{z\in N; \ (n_1-n_1^1)^2/\delta^2+\cdots+(n_k-n_k^1)^2
\leq r^2\}$ increase and first touch $b\Phi_{nr}$ for $\delta_{sup}>1$
at points $z_0\in N$ with $n_1(z_0) \neq n_1^1$. Therefore,
$B_{\delta_{sup}}$ contains open convex cones $\Gamma$ at
$(z_0,-X(z_0))$ and since $\hbox{Int} \ B_{\delta_{sup}} \subset
N\backslash \Phi_{nr}$, $\Gamma \cap \Phi_{nr}=\emptyset$.

We shall show that such a point $z_0$ is removable, thus deriving a
contradiction, firstly using the minimalization theorem and secondly
using the deformation of discs technique.

According to the Deformation Lemma, or Proposition 5.5 below, we can
assume that we are given a continuous function $f\in
\hbox{CR}(M\backslash \Phi_{nr})$ which extends holomorphically into a
wedge ${\cal W}_0$ attached to $M\backslash \Phi_{nr}$.

There exists a truncated open convex cone $C$ in $\C^n$ with vertex
$z_0$ such that $C \cap N = \Gamma$. If $N$ is $C^{2, \beta}$, $\beta
>0$, minimalize $N$ into $C$ in a manifold $M_1$ as in the proof of
Theorem 5.A and apply the argument of sweeping out by wedges. Then
$b\Phi_{nr} \backslash \{z_0\}$ is removable near $z_0$, and therefore
also $z_0$ by Theorem 1, a contradiction.

Second, assume $N$ is $C^2$-smooth.  Since there exists a cone $\Gamma
\subset N$ at $(z_0, -X(z_0))$ and $X(z_0)\in T_{z_0}^cN$, there
exists a disc $A$ attached to $N$ with $A(1) \in N$ and $A(\gamma)
\subset \Gamma$ for some open arc $\gamma \subset b \Delta$ with $-1
\in \gamma$. Indeed, choosing a holomorphic coordinate system $(w,z)$
at $z_0$ as in $(10)$ with $z_0=0$ and $T_{z_0} N=\{y=0, v_1=0\}$ and
$X(0)$ directed in the positive $u_2$-axis, this is true for the disc
with $(w_2,...,w_p)$ $\zeta$-holomorphic component equal
$(c(1-\zeta),0,...,0)$, $c>0$ small and $(w_1,z)(\zeta)$ satisfying
Bishop's equation relative to $N$.

One can introduce a manifold $K\subset N$ with $0\in K,
\frac{d}{d\theta}|_{\theta=0} A(e^{i\theta}) \not\in T_0K$,
$\hbox{codim}_N K =1$ and deformations $A_t$ of $A$ in the normal
space to $N$ at $A(-1)$ depending on a real parameter $t\in {\cal
T}\subset \R^q$, $q=\hbox{codim} \ N -1$, and on $p_0 \in {\cal K}
\subset K$, such that each $A_{t,p_0}$ is attached to $N \cup {\cal
W}_0$ and ${\cal W}_A=\{A_{t, p_0}(\stackrel{\circ}{\Delta}_1); \ t
\in {\cal T}, p_0 \in {\cal K} \}$ generates a wedge of edge $N$ at
$0$. Include $N$ in a $C^{2, \alpha}$-smooth one-parameter family
$N_s, |s|<\delta, \delta >0$ of $C^{2, \alpha}$ generic manifolds
$N_s$ contained in $M$ such that $N_0=N, N_s \cap \Phi_{nr} =
\emptyset$ for $s\neq 0$,$\cup_s N_s= M$ near $0$ and the $N_s$
contain some fixed point in $\Gamma$. Then, according to the smooth
dependence of the solutions of Bishop's equation on parameters, for
each small $s\neq 0$, the sets ${\cal W}_s
=\{A_{s,t,p_0}(\stackrel{\circ}{\Delta}_1); t\in {\cal T}, p_0 \in
{\cal K}\}$ are wedges of edge $N_s$ to which $f|_{N_s}$, which is CR
on $N_s$, holomorphically extend and whose direction depends smoothly
on $s$. As in the proof of the sweeping out by wedges lemma, the union
of the ${\cal W}_s$ for $s\neq 0$ generates a wedge ${\cal W}$ of edge
$M$ at $z_0$ and the extensions obtained stick in a well-defined
holomorphic function into ${\cal W}$.

The proof of Theorem 5.4. is complete.

\smallskip
{\it Remark.} The second proof of Theorem 5.4 could also provide a
second proof of Theorem 5.A. with $N$ of class $C^2$.

\bigskip
{\it Proof of Theorem 4.} The main observation is resumed in
Proposition 5.5 below, according to which all notions of removability
considered during the course are in fact one and the same.

Given two wedges ${\cal W}_1$ and ${\cal W}_2$, one sets ${\cal W}_1
\subset \subset {\cal W}_2$ if their cones satisfy $C_1 \cap S^{2n-1}
\subset \subset C_2 \cap S^{2n-1}$, where $S^{2n-1}$ denotes the unit
sphere in $\C^n$ identified with $R^{2n}$.  Given two wedges ${\cal
W}_0$ and ${\cal W}_0'$ attached to $M$, one sets ${\cal W}_0 \subset
\subset {\cal W}_0'$ if ${\cal W}_{0,z} \subset \subset {\cal W}_{0,
z'}$ for each $z$. Notice that one has ${\cal W}_0 \subset \subset
{\cal W}_0'$ provided ${\cal W}_0 \cap {\cal W}_0'$ contains a
nonempty wedge attached to $M$, even if ${\cal W}_0 \not\subset {\cal
W}_0'$ This is because a wedge attached to $M$ is not supposed to be
exactly a wedge of edge $M$ at each point of $M$.

\smallskip
\noindent
{\sc Proposition 5.5.} {\it Let $M$ be generic, $C^{2,
\alpha}$-smooth, let $\Phi$ be a proper closed subset of $M$, $z_0 \in
b \Phi$ and assume that $M$ is minimal at $z_0$. Then the following
are equivalent.

\begin{minipage}[t]{16cm}
\noindent
$(i)$ Given a wedge ${\cal W}_0$ attached to $M$, $\exists U \ni z_0$,
$\exists {\cal W}$ attached to $U$ such that $\forall f \in {\cal
O}({\cal W}_0),$ $\exists F \in {\cal O}({\cal W})$ with $F=f$ into a
subwedge ${\cal W}_1 \subset \subset {\cal W}_0 \cap {\cal W}$
attached to $U\backslash \Phi$.

\noindent
$(ii)$ Given a wedge ${\cal W}_0$ attached to $M\backslash \Phi$,
$\exists U\ni z_0, \exists {\cal W}={\cal W}(U, z_0)$, such that
$\forall f \in {\cal O}({\cal W}_0)$ continuous up to $M\backslash
\Phi$, $\exists F \in {\cal O}({\cal W})$ continuous up to
$U\backslash \Phi$ with $F=f$ in $U\backslash \Phi$.

\noindent
$(iii)$ Given a wedge ${\cal W}_0$ attached to $M\backslash \Phi$,
$\exists U \ni z_0$, $\forall \ V$ neighborhood of $\Phi \cap U$ in
$U$, $\forall \varepsilon >0$, $\exists U^d$ $C^{2, \alpha}$
deformation of $U$, $U^d \cap (U\backslash V)=U\backslash V$,
$||U^d-U||_{C^{2, \alpha}} <\varepsilon$ such that $\forall f \in
{\cal O}({\cal W}_0)$ continuous up to $M\backslash \Phi$, $\exists
f^d \in \hbox{CR}(U^d)$ with $f^d=f$ on $U\backslash V$.
\end{minipage}

\vspace{0.16cm}
\noindent
When these properties hold, there exists a neighborhood ${\cal U}$ of
$z_0$ in $\C^n$ and ${\cal W}$ a wedge attached to
$(M\backslash\Phi)\cup ({\cal U} \cap M)$ such that ${\cal W}\equiv
{\cal W}_0$ outside ${\cal U}$, ${\cal W}\cap {\cal W}_0$ contains a
wedge ${\cal W}_1$ attached to $M\backslash \Phi$ and for each
function $f\in {\cal O}({\cal W}_0)$ there exists $F\in {\cal O}({\cal
W})$ with $F=f$ into ${\cal W}_1$.  Moreover, these properties are
equivalent to $(i)'$, $(ii)'$, $(iii)'$ where one replaces ${\cal
W}_0$ with a neighborhood $\omega$ of $M\backslash \Phi$ in $\C^n$.}

\medskip
\noindent
{\it Proof.}  $(i) \Rightarrow (ii)$. Let ${\cal W}_0$ be wedge
attached to $M\backslash \Phi$. By $(i)$, there exists a wedge ${\cal
W}={\cal W}(U',z_0)$ such that $\forall f \in {\cal O}({\cal W}_0)$,
$\exists F \in {\cal O}({\cal W})$, $F=f$ into a subwedge ${\cal W}_1
\subset \subset {\cal W}_0 \cap {\cal W}$ attached to $U\backslash
\Phi$.  Since ${\cal W}_1 \subset \subset {\cal W}_0$, for such $f$
continous up to $M\backslash \Phi$, $f|_{{\cal W}_1}$ admits a
continuous limit up to $U\backslash \Phi$, and then also $F|_{{\cal
W}_1}=f|_{{\cal W}_1}$, proving that $F$ is continous up to
$U\backslash \Phi$.  Hence $(ii)$ holds.

$(ii) \Rightarrow (iii)$. Let ${\cal W}_0$ be a wedge attached to
$M\backslash \Phi$.  By $(ii)$, there exists a wedge ${\cal W}={\cal
W}(U', z_0)$ over $U'$ at $z_0$ such that $\forall f \in {\cal
O}({\cal W}_0)$ continous up to $M\backslash \Phi$, $\exists F \in
{\cal O}({\cal W})$ continuous up to $U'\backslash \Phi$ with $F=f$ in
$U' \backslash \Phi$. Let $U \subset \subset U'$ be a subneighborhood,
let $V$ be a neighborhood of $\Phi \cap U$ in $U$ and $\varepsilon >0$
arbitrary. Since ${\cal W}(U', z_0)$ is a wedge of edge $U'$ at $z_0$,
there exsists a $C^{2, \alpha}$ deformation $U^d\subset {\cal W} \cup
(U\backslash V)$ of $U$ with $U^d\equiv U$ in $U\backslash V$ and
$||U^d - U||_{C^{2, \alpha}}<\varepsilon$. For each function $f$
extending as a $F \in {\cal O}({\cal W})$, one sets $f^d:= F|_{U^d}$
which is CR without singularities on $U^d$ and satisfies $f^d=f$ on
$U\backslash V$.  Hence $(iii)$ holds.

$(iii) \Rightarrow (i)$. Let ${\cal W}_0$ be a wedge attached to
$M\backslash \Phi$.  Choose a wedge ${\cal W}_1 \subset {\cal W}_0$
attached to $M\backslash\Phi$ with ${\cal W}_1\subset \subset {\cal
W}_0$ near $z_0$, a $C^{2, \alpha}$ varying direction on $M\backslash
\Phi$ near $z_0$ and ${\cal W}_1 \equiv {\cal W}_0$ away from $z_0$.
Deform first $M$ in a one-parameter family of $C^{2, \alpha}$
manifolds $M^{d_1}, d_1 \geq 0, M^{d_1} \subset {\cal W}_1 \cup \Phi$
with $||M^{d_1}-M||_{C^{2, \alpha}} <\varepsilon(d_1) < \delta/2,
\delta>0, \varepsilon(d_1) \to 0$ as $d_1 \to 0$ and for each function
$f\in {\cal O}({\cal W}_0)$ take $f^{d_1}=f|_{M^{d_1}}\in
\hbox{CR}(M^{d_1}\backslash \Phi)$. Then the $f^{d_1}$ are in fact
defined in a neighborhood ${\cal W}_1$ of $M^{d_1}\backslash \Phi$ in
$\C^n$ and holomorphic there.  Choose $\delta>0$ so small that for
each $C^{2, \alpha}$ manifold $M'$ with $||M'-M||_{C^{2,
\alpha}}<\delta$, there exists a wedge ${\cal W}(A')$ associated with
a disc of zero defect $A'$ attached to $M'$ which is a smooth
perturbation of a disc $A$ attached to $M$ through $z_0$ of zero
defect and the size $U'$ of the base of ${\cal W}(A')$ in $M'$
satisfies $\hbox{dist}(bU', z_0') \geq \kappa >0$ for every $M'$.

According to $(iii)$, there exists a neighborhood $U^{d_1}\ni z_0$
such that, given $V$ a neighborhood of $\Phi\cap U^{d_1}$ in $U^{d_1}$
and $\varepsilon_2 >0$ arbitrary, there exists a $C^{2, \alpha}$
deformation $(U^{d_1})^{d_2}$ of $U^{d_1}$, $d_2\geq 0,
(U^{d_1})^0=U^{d_1}$ with $||(U^{d_1})^{d_2}-U^{d_1}||_{C^{2,
\alpha}}<\varepsilon_2 <\delta/2$ such that $\forall f \in {\cal
O}({\cal W}_1), \exists f^{d_2} \in \hbox{CR}((U^{d_1})^{d_2})$ with
$f^{d_2}\equiv f^{d_1}$ in $U^{d_1}\backslash
V=(U^{d_1})^{d_2}\backslash V$.  Since $f^{d_2}$ is CR and we can
assume that the size of $A^{d_1}$ is smaller than $U^{d_1}$, so there
exist perturbations $(A^{d_1})^{d_2}$ of $A$ of zero defect attached
to $(U^{d_1})^{d_2}$, therefore $f^{d_2}$ extends holomorphically into
a wedge ${\cal W}((A^{d_1})^{d_2})$ of edge $(U^{d_1})^{d_2}$ near
$(z_0^{d_1})^{d_2}$.  Letting $d_2, \varepsilon_2$ tend to zero and
$V$ shrink to $\Phi$, this shows that $f^{d_1}\in
\hbox{CR}(U^{d_1}\backslash\Phi)$ extends holomorphically into ${\cal
W}(A^{d_1})$ as a function $F^{d_1} \in {\cal O}({\cal W}(A^{d_1}))$
continuous up to $U^{d_1}\backslash \Phi$.  In fact, we just showed
that $(iii)\Rightarrow (ii)$.

Then the $F^{d_1}$ stick together in a well-defined holomorphic
function in ${\cal W}_2= \bigcup_{d_1>0} {\cal W}(A^{d_1})$. Indeed,
the wedges ${\cal W}(A^{d_1})$ varying differentiably as
$\frac{\partial A^{d_1}}{\partial \zeta}(1)$, they make successive
connected intersection and one obtains a function $F^{d_2} \in {\cal
O}({\cal W}_2)$.  Shrinking the height of ${\cal W}_1$, {\em i.e.} the
troncature of the cones defining ${\cal W}_1$ at each point, one can
insure that ${\cal W}_2 \cap {\cal W}_1$ has as many connected
components as $U\backslash \Phi$, for some neighborhood $U\ni z_0$ in
$M$.  Thus, $f|_{{\cal W}_1}$ and $F_2$ stick in a well-defined
function $F\in {\cal O}({\cal W}), {\cal W}:={\cal W}_1\cup {\cal
W}_2$. Hence $(i)$ holds.

The functions $f\in {\cal O}({\cal W}_1)$ and $F\in {\cal O}({\cal
W})$ stick together in a holomorphic function in the wedge ${\cal W}_1
\cup {\cal W}$ attached to $(M\backslash \Phi)\cup U_1$, since
$f\equiv F$ in each connected component of ${\cal W}_1$.

The equivalence of $(i)', (ii)'$ and $(iii)'$ is proved in exactly the
same way. $(iii)$ implies $(iii)'$ and conversely, in the proof of
$(iii)\Rightarrow (i)$, it was in fact proved that $(iii)'\Rightarrow
(i)$.

The proof of Proposition 5.5 is complete.

We can also translate the equivalences of Proposition 5.5 when one
assumes that only a CR function is previously given.

\smallskip
\noindent
{\sc Corollary 5.6.} {\it Let $M$ be generic, $C^{2,\alpha}$-smooth,
let $\Phi\subset M$ be a proper closed subset of $M$, $z_0\in b \Phi$
and assume that $M$ is minimal at $z_0$. Then the following are
equivalent.

\begin{minipage}[t]{16cm}
\noindent
$(iv)$ $\exists$ ${\cal W}$ of edge $M$ at $z_0$ such that $\forall f
\in \hbox{CR}(M\backslash \Phi)$ $\exists F \in {\cal O}({\cal W})$
continuous up to $M\backslash \Phi$ with $F=f$ there.

\noindent
$(v)$ $\exists U \ni z_0$ such that $\forall V$ neighborhood of $\Phi
\cap U$ in $U$, $\forall \varepsilon >0$, $\exists U^d$ $C^{2,
\alpha}$, $U^d \cap (U\backslash V)=U\backslash V$, $||U^d-U||_{C^{2,
\alpha}}<\varepsilon$ such that $\forall f \in \hbox{CR}(M\backslash
\Phi)$, $\exists f^d\in \hbox{CR}(U^d)$, $f^d=f$ on $U\backslash V$.
\end{minipage}

\vspace{0.16cm}
\noindent
If $\hbox{CR}(M\backslash \Phi)$ has the wedge extension property at
every point of $M\backslash \Phi$, these are still equivalent to $(i),
(ii), (iii)$ above.}

\smallskip
{\it Remark.}  Assume now that $\Phi=N$ is a proper at least
one-codimensional submanifold of $M$. We can avoid assuming that $M$
is minimal at $z_0$ in the equivalence above. Using normal
deformations of discs as in the proof of Theorem 5.4, we could prove
also:

\smallskip
\noindent
{\sc Corollary 5.7.} {\it Let $M$ be generic, $C^{2, \alpha}$-smooth,
let $N\subset M$ be a submanifold, $\hbox{codim}_M N \geq 1$, let
$z_0\in N$ and assume that $T_{z_0} N \not\supset T_{z_0}^c M$. Then
the equivalences of Proposition 5.5 hold with $\Phi=N$.}

\smallskip
{\it Proof.} We only suggest that, when assuming that $T_{z_0} N
\not\supset T_{z_0}^cM$, we implicitely assume that there exists a
disc $A$ attached to $M$ with $A(1)=z_0\in N$ and $A(\gamma)\cap N =
\emptyset$ for some open arc $\gamma\ni -1$ in $b\Delta$. Then,
$(iii)$ enables one to consider CR functions on a deformation of $M$
without singularities and the deformations of $A$ in the normal space
to $M$ at $A(-1)$ in ${\cal W}_0$ render it possible that $A$ plays
the role of a disc of zero defect as in the proof of $5.5,
(iii)\Rightarrow (i)$.

The proof of Corollary 5.7 is complete.

\bigskip
{\it Proof of Theorem 5.} Since $T_{z_0}N + T_{z_0}^cM=T_{z_0}M$ and
$\hbox{codim}_M N =2$, we can assume that $N$ is given in a coordinate
system as $(10)$ by
$$y=h(w,x) \ \ \ \ \ w_p=g(w',x),$$ where $w'=(w_1,...,w_{p-1})$,
$h(0)=0, dh(0)=0$ and $g(0)=0$. Since $N$ is CR and has $\hbox{CRdim}
\ N=p-1$, $g$ is CR on the manifold $N^{\pi}\subset\C^{p-1}$ with
equation $y=H(w',x)$, where $H(w',x)=h(w',g(w',x),x)$.  Since
$N^{\pi}$ is minimal at $0$, there exists a family $A_{t, n_0}^{\pi}$
of discs attached to $N^{\pi}$ in $\C^{n-1}$ such that for $t$ in a
neighborhood ${\cal T}$ of $0$ in $\R^{q-1}$, $n_0$ in a one
codimensional manifold ${\cal K}\subset N^{\pi}$ through $0$ and
$\zeta\in \stackrel{\circ}{\Delta}_1$, the set of points
$A_{t,n_0}^{\pi}(\zeta)$ spans a wedge ${\cal W}^{\pi}$ of edge
$N^{\pi}$ in $\C^{n-1}$. We can assume that the $A_{t, n_0}^{\pi}$
have a $w'$-component which embeds $\overline{\Delta}$ into $\C^{n-1}$
and $\frac{d}{d\theta}|_{\theta=0} A(e^{i\theta}) \not\in
T_0^cN^{\pi}$, since $N^{\pi}$ is minimal at every point (see the
argument in the proof of Proposition 2.1). Set
$A^{\pi}_{t,n_0}(\zeta)=(w_{t, n_0}'(\zeta), z_{t,n_0}(\zeta))$.
Since $g$ is CR, $A_{t,n_0}(\zeta)=(w_{t,n_0}'(\zeta),
g(w_{t,n_0}'(\zeta), x_{t,n_0}(\zeta)), z_{t, n_0}(\zeta))$ is a
holomorphic disc in $\zeta$ attached to $N$.  For $a$ in a
neighborhood ${\cal V}$ of $0$ in $\C$, consider the analytic disc
attached to $M$ $A_{t,n_0,a}(\zeta)=(w_{t,n_0}'(\zeta),
g(w_{t,n_0}'(\zeta),x_{t,n_0}(\zeta))+a, z_{t,n_0,a}(\zeta))$ where
$x_{t,n_0,a}=-T_1h(w_{t,n_0}',
g(w_{t,n_0}',x_{t,n_0})+a,x_{t,n_0,a})+x_{t,n_0}(1)$.  Since $h(0)=0,
dh(0)=0$, we have the estimate $|x_{t,n_0,a}-x_{t,n_0}|< \varepsilon
|a|$, $\varepsilon <<1$ and this proves that $A_{t,n_0,a}(b\Delta)\cap
N= \emptyset$ when $a\neq 0$. Furthermore, if the size of
$A_{t,n_0}^{\pi}$ is sufficiently small, by taking $a=a(s)\neq 0$,
$0\leq s \leq 1$, one can check that the discs $A_{t, n_0,a}$ are
analytically isotopic to a point in $M\backslash N$. Since $\hbox{dim}
\ T_{z_0}^cM/ (T_{z_0}^cM \cap T_{z_0}N)=2$, $M\backslash N$ is
globally minimal around $z_0$, hence CR functions on $M\backslash N$
are wedge extendible at every point of $M\backslash N$ near $z_0$
\cite{ME}.  Therefore the argument in the proof of Proposition 4.2 can
be repeated here: CR functions are holomorphically extendible into the
wedge open set
\begin{equation} {\cal W}\backslash {\cal W}^{an}=\{A_{t,n_0,a}(\zeta);\ 
t\in {\cal T}, n_0 \in {\cal K}, a \in {\cal V}\backslash \{0\}, \zeta
\in \stackrel{\circ}{\Delta}_1 \}
\end{equation}
minus the analytic wedge
$${\cal W}^{an}=\{A_{t,n_0}(\zeta); \ t\in {\cal T}, n_0\in {\cal K},
\zeta\in \stackrel{\circ}{\Delta}_1\}.$$

The proof of Theorem 5 is complete.

\bigskip
\begin{flushright}
{\footnotesize \it D\'epartement de Math\'ematiques et d'Informatique}
\end{flushright}
\vspace{-22pt}
\begin{flushright}
{\footnotesize \it \'Ecole Normale Sup\'erieure, 45 rue d'Ulm, F-75230
Paris Cedex 05.}
\end{flushright}
\vspace{-22pt}
\begin{flushright}
{\footnotesize \it E-mail address: {\rm merker@dmi.ens.fr}}
\end{flushright}
\end{document}